\documentclass[aos,citesort,dvips]{arximspdf}
\usepackage{dcolumn,multirow}
\usepackage{graphicx}

\doi{10.1214/09-AOS749}
\volume{38}
\issue{5}
\pubyear{2010}
\firstpage{3093}
\lastpage{3128}

\makeatletter
\newtheorem{theorem}{Theorem}
\newtheorem{assumption}{Assumption}
\newcommand{\ds}{\displaystyle}
\newproclaim{remark}{Remark}
\newcolumntype{d}[1]{D{.}{.}{#1}}
\newtheorem{lemma}{Lemma}
\newtheorem{proposition}{Proposition}
\makeatother

\begin{document}
\begin{frontmatter}

\title{Is Brownian motion necessary to model high-frequency data?}
\runtitle{Brownian motion or not}

\begin{aug}
\author[A]{\fnms{Yacine} \snm{A\"{\i}t-Sahalia}\ead[label=e1]{yacine@princeton.edu}\thanksref{t1}\corref{}}
and
\author[B]{\fnms{Jean} \snm{Jacod}\ead[label=e2]{jean.jacod@upmc.fr}}
\thankstext{t1}{Supported in part by NSF Grants DMS-05-32370 and SES-0850533.}
\runauthor{Y. A\"{\i}t-Sahalia and J. Jacod}
\pdfauthor{Yacine Ait-Sahalia, Jean Jacod}
\affiliation{Princeton University and UPMC (Universit\'e Paris-6)}
\address[A]{Department of Economics\\
Princeton University and NBER\\
Princeton, New Jersey 08544-1021\\
USA\\
\printead{e1}} 
\address[B]{Institut de Math\'ematiques de Jussieu\\
CNRS UMR 7586\\
Universit\'e P. et M. Curie (Paris-6)\\
75252 Paris C\'edex 05\\
France\\
\printead{e2}}
\end{aug}

\received{\smonth{3} \syear{2009}}
\revised{\smonth{9} \syear{2009}}

\begin{abstract}
This paper considers the problem of testing for the presence of a continuous part in
a semimartingale sampled at high frequency. We provide two tests, one where the null hypothesis is that a continuous
component is present, the other where the continuous component is absent, and the model is then driven by a pure
jump process. When applied to high-frequency individual stock data, both tests point toward the need to
include a continuous component in the model.
\end{abstract}


\end{frontmatter}

\section{Introduction}\label{sec:intro}

This paper continues our development of statistical methods designed to assess
the specification of continuous-time models sampled at high frequency. The
basic framework, inherited from theoretical models in mathematical finance but
also common in other fields such as physics or biology, is one where the
variable of interest $X$, in financial examples often the log of an asset
price, is assumed to follow an It\^{o} semimartingale. That semimartingale is
observed on some fixed time interval $[0,T]$ at discrete regularly spaced
times $i\Delta_{n}$, with a time lag $\Delta_{n}$ which is small.

A semimartingale can be decomposed into the sum of a drift, a continuous
Brownian-driven part and a discontinuous, or jump, part. The jump part can in
turn be decomposed into a sum of ``small
jumps''\ and ``big jumps.''%
\ Such a process will always generate a finite number of big jumps, but it may
give rise to either a finite or infinite number of small jumps, corresponding
to the finite and infinite jump activity situations, respectively. In earlier
work, we developed tests to determine on the basis of the observed sampled
path on $[0,T]$ whether a jump part was present, whether the jumps had finite
or infinite activity, and in the latter situation proposed a definition and an
estimator of a ``degree of jump activity'' parameter.

In this paper, we tackle the last remaining question: Does the semimartingale
need to have a continuous part? In other words, is the Brownian motion present
at all? From a model specification standpoint, there is a natural statistical
interest in distinguishing the two situations where a continuous part is
included or not, on the basis of an observed sample path. When there are no
jumps, or finitely many jumps, and no Brownian motion, $X$ reduces to a pure
drift plus occasional jumps, and such a model is fairly unrealistic in the
context of most financial data series, although it may be realistic in some
other contexts. But for financial applications one can certainly consider
models that consist only of a jump component, plus perhaps a drift, if that
jump component is allowed to be infinitely active.

Many models in mathematical finance do not include jumps. But among those that
do, the framework most often adopted consists of a jump-diffusion: these
models include a drift term, a Brownian-driven continuous part and a finite
activity jump part\ (see, e.g., \cite{balltorous83,bates91} and \cite{merton76}). When infinitely many jumps are included, however, there are a
number of models in the literature which dispense with the Brownian motion
altogether. The log-price process is then a purely discontinuous L\'{e}vy
process with infinite activity jumps or, more generally, is driven by such a
process (see, e.g., \cite{eberleinkeller95,carrwu03b} and \cite{madanseneta90}).

The mathematical treatment of models relying on pure jump processes is quite
different from the treatment of models where a Brownian motion is present. For
instance, risk management procedures, derivative pricing and portfolio
optimization are all significantly altered, so there is interest from the
mathematical finance side in finding out which model is more likely to have
generated the data.

For all these reasons, it is of importance to construct procedures which allow us
to decide whether the Brownian motion is really here, or if it can be forgone  in favor of a pure jump process. This is the aim of this paper: we will
provide two tests allowing for a symmetric treatment of the two situations
where the null hypothesis is that the Brownian motion is present, and where
the null is that the Brownian motion is absent.

In the context of a specific parametric model, allowing for jump components of
finite or infinite activity on top of a Brownian component,
\cite{carrgemanmadanyor02} find that the time series of index returns are
likely to be devoid of a continuous component. An alternative but related
approach to testing for the presence of a Brownian motion component to the one
we propose here is due to \cite{tauchentodorov}. They employ the test
statistic for jumps of \cite{yacjacod09a}, plot its logarithm for different
values of the power argument and contrast the behavior of the plot above two
and below two in order to identify the presence of a Brownian component. A
formal test is constructed under the null hypothesis where a continuous
component is present.

The methodology that both \cite{tauchentodorov} and we employ to design our
respective test statistics is based on tried-and-true principles that
originate in our earlier work on testing whether jumps are present
\cite{yacjacod09a}, whether they have finite or infinite activity
\cite{yacjacod08c} and on estimating the index of jump activity
\cite{yacjacod09b}, although, of course, exploited in a manner specific to the
problem at hand. We compute power variations of the increments, suitably
truncated and/or sampled at different frequencies. Exploiting the different
asymptotic behavior of the variations as we vary these parameters gives us
enough flexibility to accomplish our objectives. As is well known, powers
below two will emphasize the continuous component of the underlying sampled
process. Powers above two will conversely accentuate its jump component. The
power two puts them on an equal footing. Truncating the large increments at a
suitably selected cutoff level can eliminate the big jumps when needed, as was
shown by \cite{mancini01}. Finally, sampling at different frequencies can let
us distinguish between situations where the variations converge to a finite
limit, in which case the ratio of two variation measures constructed at
different frequencies will converge to one, from situations where the
variations converge to either zero or diverge to infinity, in which case the
ratio will typically converge to a different constant. Since these various
limiting behaviors are indicative of which component of the model dominates at
a particular power, they effectively allow us to distinguish between all
manners of null and alternative hypotheses.

This said, the commonality of approach should not mask the fact that each
situation is, in reality, mathematically quite different. By nature, certain
components of the model are turned off under particular null hypotheses. For
instance, when the null hypothesis is that no Brownian motion is present, as
will be the case for our first test here, then jumps drive the asymptotics. As
a result, the driving component of the model that matters for the asymptotic
behavior of the statistic will vary with the situation and consequently the
methods employed behind the scenes to obtain the desired asymptotics will vary accordingly.

The paper is organized as follows. Section \ref{sec:model} describes our model
and the statistical problem. Our testing procedure is described in Section
\ref{sec:tests}, and the next two Sections, \ref{sec:MC} and \ref{sec:data}, are
devoted to a simulation study of the tests and an empirical implementation of
our tests on high-frequency stock returns. Section \ref{sec:TR} is devoted to
technical results and to the proof of the main theorems.

\section{The model}\label{sec:model}

The underlying process $X$ which we observe at discrete times is a
$1$-dimensional It\^{o} semimartingale defined on some filtered space
$(\Omega,\mathcal{F},(\mathcal{F}_{t})_{t\geq0},\mathbb{P})$, which means that
its characteristics $(B,C,\nu)$ are absolutely continuous with respect to
Lebesgue measure. $B$ is the drift, $C$ is the quadratic variation of the
continuous martingale part and $\nu$ is the compensator of the jump measure
$\mu$ of $X$. In other words, we have
\begin{eqnarray}\label{9}
B_{t}(\omega)&=&\int_{0}^{t}b_{s}(\omega)\,ds,\qquad
C_{t}(\omega)=\int_{0}^{t}\sigma_{s}(\omega)^{2}\,ds,\nonumber\\[-8pt]\\[-8pt]
\nu(\omega,dt,dx)&=&dt \,F_{t}(\omega,dx).\nonumber
\end{eqnarray}
Here $b$ and $\sigma$ are optional process, and $F=F_{t}(\omega,dx)$ is a
transition measure from $\Omega\times\mathbb{R}_{+}$ endowed with the
predictable $\sigma$-field into $\mathbb{R}\setminus\{0\}$. More customarily,
one may write $X$ as
\begin{eqnarray}\label{1}
X_{t}  &
=&X_{0}+\int_{0}^{t}b_{s}\,ds+\int_{0}^{t}\sigma_{s}\,dW_{s}\nonumber\\[-8pt]\\[-8pt]
&&{}  +\int_{0}^{t}\int x1_{\{|x|\leq1\}}(\mu-\nu)(ds,dx)+\int_{0}^{t}\int
x1_{\{|x|>1\}}\mu(ds,dx),\nonumber
\end{eqnarray}
where $W$ is a standard Brownian motion. It is also possible to write the last
two terms above as integrals with respect to a Poisson measure and its
compensator, but we do not need this here. This is a standard setup and we
refer the reader to \cite{jacodshiryaev2003} for details.

We have referred above to ``small jumps''\ and
``big jumps.''\ In the context of (\ref{1}),
they are represented, respectively, by the last two integrals. The size cutoff
$1$ adopted here is arbitrary and could be replaced by any fixed
$\varepsilon>0,$ a change which amounts merely to an adjustment to the drift
term $B_{t}.$ Note that the small jumps integral needs to be compensated by
$\nu$ since there are potentially an infinite number of such small jumps. The
large jump integral is always a finite sum; it may be compensated if desired
but this is not necessary. Any compensation or lack thereof is then again
absorbed by an adjustment to the drift.

We now turn to the assumptions. As usual for tests, the assumptions
essentially ensure that one can compute and then estimate a significance level
under the null hypothesis. So here, we need some structure for the jumps of
$X$, namely that the small jumps essentially behave like the small jumps of a
stable process with some index $\beta$, up to a random intensity. As noted
above, when no Brownian is present, we view the realistic situation as one
where there are infinitely many small jumps. When the null is that there is a
Brownian motion, we need the additional assumption that the volatility process
$\sigma_{t}$ is itself an It\^o semimartingale.

We would like to give tests with a prescribed asymptotic level, as
$n\rightarrow\infty$, and, of course, this is more difficult when $\beta$
increases because then the process resembles more and more a continuous
process plus a few big jumps: The qualitative behavior of the paths can become
quite similar whether the Brownian motion is present or not. So,
unsurprisingly, we can exhibit a test with prescribed level, for the null
hypothesis where the Brownian motion is present, only when $\beta<1$. The
parameter $\beta$ is typically unknown (although a method for estimating
$\beta$ in this setting is given in \cite{yacjacod09b}). On the other hand,
for the null hypothesis where the Brownian motion is absent\ we provide a test
which works under no assumption on $\beta$.

With this context in mind, here is the first assumption which will be assumed throughout:

\begin{assumption}\label{ass:A1}
\textup{(i)} The drift process $b_{t}$ is locally bounded and the
volatility process $\sigma_{t}$ is c\`{a}dl\`{a}g.

\textup{(ii)} There are three constants $0\leq\beta^{\prime\prime}\leq\beta^{\prime
}<\beta<2$ and a locally bounded process $L_{t}\geq1$, such that the L\'{e}vy
measure $F_{t}$ is of the form $F_{t}=F_{t}^{\prime}+F_{t}^{\prime\prime}$,
where
\begin{equation} \label{5}
\quad F_{t}^{\prime}(dx) = \frac{\beta(1+|x|^{\beta-\beta^{\prime}}f(t,x))}%
{|x|^{1+\beta}} \bigl(  a_{t}^{(+)}1_{\{0<x\leq z_{t}^{(+)}\}}+a_{t}%
^{(-)}1_{\{-z_{t}^{(-)}\leq x<0\}} \bigr)\,  dx,
\end{equation}
where $a_{t}^{+}$, $a_{t}^{-}$, $z_{t}^{+}$, $z_{t}^{-}$ are nonnegative
predictable processes and $f=f(\omega,t,x)$ is predictable function (meaning
$\mathcal{P}\otimes\mathcal{B}(\mathbb{R})$-measurable, where $\mathcal{P}$ is
the predictable $\sigma$-field on $[0,\infty)\times\Omega$), satisfying
\begin{eqnarray}\label{BA2}
\frac{1}{L_{t}}&\hspace*{3pt}\leq& z_{t}^{(+)}\leq1,\qquad
\frac{1}{L_{t}}\leq z_{t}^{(-)}\leq1,\nonumber\\[-8pt]\\[-8pt]
A_{t}&:=&a_{t}^{(+)}+a_{t}^{(-)}\leq L_{t},\qquad
|f(t,x)|\leq L_{t},\nonumber
\end{eqnarray}
and where $F_{t}^{\prime\prime}$ is a measure which is singular with respect
to $F_{t}^{\prime}$ and satisfies
\begin{equation}\label{7}
\int(|x|^{\beta^{\prime\prime}}\wedge1) F_{t}^{\prime\prime}(dx)\leq L_{t}.
\end{equation}
\end{assumption}

This assumption is identical to Assumptions 1 and 2 of \cite{yacjacod09b}
[with some notational changes: $(\gamma,\beta^{\prime},a_{t}^{+},a_{t}^{-})$
in that paper are called here $(\beta-\beta^{\prime},\beta^{\prime\prime
},\beta a_{t}^{+},\break\beta a_{t}^{-})$, and the condition $\beta^{\prime\prime
}\leq\beta^{\prime}$ is not a restriction and is put here only for convenience].

For example, take a process solution of the stochastic differential equation%
\begin{equation}\label{eq:exmodel}
dX_{t}=b_{t}\,dt+\sigma_{t}\,dW_{t}+\delta_{t-}\,dY_{t}+\delta_{t-}^{\prime}%
\,dY_{t}^{\prime},
\end{equation}
where $\delta$ and $\delta^{\prime}$ are c\`{a}dl\`{a}g adapted processes, $Y$
is $\beta$-stable or tempered $\beta$-stable and $Y^{\prime}$ is any other
L\'{e}vy process whose L\'{e}vy measure integrates $|x|^{\beta^{\prime\prime}%
}$ near the origin and has an absolutely continuous part whose density is
smaller than $K|x|^{-(1+\beta^{\prime})}$ on $[-1,1]$ for some $K>0$ (e.g., a stable process with index strictly smaller than $\beta^{\prime}$).
Then $X$ will satisfy Assumption \ref{ass:A1}.

If this assumption is satisfied with $\beta<1$, then almost surely the jumps
have finite variation $\sum_{s\leq t}|\Delta X_{s}|<\infty$ for all $t$ or
equivalently, $\int_{0}^{t}\int|x|\mu(ds,dx)<\infty$. This allows us to decompose
$X$ into the sum $X=X^{\prime}+X^{\prime\prime}$, where
\begin{equation}\label{701}
X_{t}^{\prime} = X_{0}+\int_{0}^{t}b_{s}^{\prime} \,ds+\int_{0}^{t}\sigma
_{s} \,dW_{s},\qquad   X_{t}^{\prime\prime} = \sum_{s\leq t}\Delta X_{s},
\end{equation}
and where $b_{t}^{\prime}=b_{t}-\int x1_{\{|x|\leq1\}} F_{t}(dx)$ is a
locally bounded process.

For clarity, we will derive the properties of both tests under the same
generic Assumption \ref{ass:A1} even though the properties of the test for the
null of a Brownian present remain valid under weaker assumptions. When the
null hypothesis to be tested is that the Brownian motion is present, it
becomes the driving process for our test statistic and as is customary for
tests or estimation problems involving a stochastic volatility, we then need
an additional regularity assumption on the $\sigma$ process:

\begin{assumption}
\label{ass:A2} We have Assumption \ref{ass:A1} with $\beta<1$. Moreover the
volatility process $\sigma_{t}$ is an It\^{o} semimartingale, that is, it can
be written (necessarily in a unique way) as
\begin{equation}
\label{2}\sigma_{t}=\sigma_{0}+\int_{0}^{t}\widetilde{b}_{s}\,ds+\int_{0}^{t}
\widetilde{\sigma}_{s}\,dW_{s}+N_{t}+\sum_{s\leq t}\Delta\sigma_{s}
1_{\{|\Delta\sigma_{s}|>1\}},
\end{equation}
where $N$ is a local martingale which is orthogonal to the Brownian motion
$W$, and further the compensator of the process $[N,N]_{t}+\sum_{s\leq
t}1_{\{|\Delta\sigma_{s}|>1\}}$ is of the form $\int_{0}^{t}n_{s}\,ds$. Moreover
we suppose that:
\begin{longlist}[(ii)]
\item[(i)] the processes $\widetilde{b}_{t}$ and $n_{t}$ are locally bounded;

\item[(ii)] the processes $\widetilde{\sigma}_{t}$ and $b^{\prime}_{t}$ defined above
are c\`{a}dl\`{a}g.
\end{longlist}
\end{assumption}

\section{The two tests}\label{sec:tests}

\subsection{The hypotheses to be tested}

In a semimartingale model like (\ref{1}), saying that the Brownian motion $W$
is absent on the interval $[0,T]$ does not mean that there is no Brownian
motion on the probability space (something which cannot be tested at all,
obviously) but it means that the Brownian motion does not impact the observed
process $X$, in the sense that the corresponding stochastic integral vanishes
on this interval, or equivalently $\sigma_{s}=0$ for Lebesgue-almost all $s$
in $[0,T]$, and it would be more appropriate to say that we are testing
whether ``the continuous martingale part of $X$ vanishes on
$[0,T]$, or not.''   This is typically an $\omega$-wise
property: we can divide the set $\Omega$ into two complementary subsets
\begin{equation}\label{6}
\Omega_{T}^{W} =  \biggl\{  \int_{0}^{T}\sigma_{s}^{2}\,ds>0 \biggr\}
,\qquad  \Omega_{T}^{\mathit{noW}} =  \biggl\{  \int_{0}^{T}\sigma_{s}^{2}\,ds=0 \biggr\}  .
\end{equation}
Then almost surely on the set $\Omega_{T}^{\mathit{noW}}$ the integral process
$X_{t}^{c}=\int_{0}^{t}\sigma_{s}\,dW_{s}$ vanishes on $[0,T]$, whereas it does
not vanish on the complement $\Omega_{T}^{W}$. In what follows, we take
$\Omega_{T}^{W}$ to represent the hypothesis that the Brownian motion is
present and $\Omega_{T}^{\mathit{noW}}$ to represent the hypothesis that the Brownian
motion is not present.

In connection with Assumption \ref{ass:A1} we consider the following set
representing paths that have infinite jump activity of some index $\beta
\in (  0,2 )  $:
\begin{equation}\label{BG6}
\Omega_{T}^{i\beta} = \{\overline{A}_{T}>0\} \qquad \mbox{where }  \overline
{A}_{t} = \int_{0}^{t}A_{s}\,ds.
\end{equation}
One knows that on the set $\Omega_{T}^{i\beta}$ the path of $X$ over $[0,T]$
has almost surely infinitely many jumps.

We are interested in testing the following two situations:
\begin{equation}\label{T-100}
\cases{
H_{0}\dvtx\Omega_{T}^{W} & vs.\quad $H_{1}\dvtx\Omega_{T}^{\mathit{noW}}$,\cr
H_{0}\dvtx\Omega_{T}^{\mathit{noW}} & vs. \quad $H_{1}\dvtx\Omega_{T}^{W}$.
}
\end{equation}
As discussed above, the realistic situation supposes that infinite activity
jumps are present when under $\Omega_{T}^{\mathit{noW}}$ and so we will in fact provide
testing procedures for the following two situations:
\begin{equation} \label{T-101}
\cases{
H_{0}\dvtx\Omega_{T}^{W} & vs. \quad $H_{1}\dvtx\Omega_{T}^{\mathit{noW}}\cap \Omega_{T}^{i\beta}$,\cr
H_{0}\dvtx\Omega_{T}^{\mathit{noW}}\cap\Omega_{T}^{i\beta} & vs. \quad
$H_{1}\dvtx\Omega_{T}^{W}\cap\Omega_{T}^{i\beta}$.
}
\end{equation}
In the second test, requiring $\Omega_{T}^{i\beta}$ with $\Omega_{T}^{W}$
under $H_{1}$ allows us to characterize precisely the properties of the
statistic under this alternative (as opposed to just $\Omega_{T}^{W}$). But it
is not necessary for the actual implementation of the test which relies on
its behavior under the null.

Finally, we recall that testing a null hypothesis ``we are in
a subset $\Omega_{0}$''\ of $\Omega$, against the alternative
``we are in a subset $\Omega_{1}$,'' with, of
course, $\Omega_{0}\cap\Omega_{1}=\varnothing$, amounts to finding a critical
(rejection) region $C_{n}\subset\Omega$ at stage $n$. The asymptotic size and
asymptotic power for this sequence $(C_{n})$ of critical regions are the
following numbers:
\begin{equation}\label{T-1}
 \cases{
\ds a = \sup \Bigl(  \limsup_{n} \mathbb{P}(C_{n}\mid A)\dvtx A\in\mathcal{F}%
, A\subset\Omega_{0}, \mathbb{P}(A)>0 \Bigr), \cr
\ds P = \inf \Bigl(  \liminf_{n} \mathbb{P}(C_{n}\mid A)\dvtx A\in\mathcal{F}%
, A\subset\Omega_{1}, \mathbb{P}(A)>0 \Bigr)  .
}
\end{equation}

\subsection{The building blocks}

Before stating the results, we introduce some notation to be used throughout.
We observe the increments of $X$
\begin{equation}\label{T-3}
\Delta_{i}^{n}X = X_{i\Delta_{n}}-X_{(i-1)\Delta_{n}},
\end{equation}
to be distinguished from the (unobservable) jumps of the process, $\Delta
X_{s}=X_{s}-X_{s-}$. In a typical application, $X$ is a log-asset price, so
$\Delta_{i}^{n}X$ is the recorded log-return over $\Delta_{n}$ units of time.

For any given cutoff level $u>0$ we count the number of increments of $X$ with
size bigger than $u$, that is,
\begin{equation}\label{703}
U(u,\Delta_{n})_{t}=\sum_{i=1}^{[t/\Delta_{n}]}1_{\{|\Delta_{i}^{n}X|>u\}}.
\end{equation}
If $p>0$ we also sum the $p$th absolute power of the increments of
$X$, truncated at level $u$, that is,
\begin{equation} \label{T-4}
B(p,u,\Delta_{n})_{t}=\sum_{i=1}^{[t/\Delta_{n}]}|\Delta_{i}^{n}%
X|^{p} 1_{\{|\Delta_{i}^{n}X|\leq u\}}.
\end{equation}
$B$ is what we call a ``truncated power
variation.'' Note that in $B$ we are retaining all increments
smaller than $u,$ whereas in $U$ we are retaining those larger than $u.$

We take a sequence $u_{n}$ of positive numbers, which will serve as our
thresholds or cutoffs for truncating the increments, and will go to $0$ as the
sampling frequency increase. There will be restrictions on the rate of
convergence of this sequence, expressed in the form
\begin{equation}\label{M-7}
u_{n}/\Delta_{n}^{\rho_{-}}\rightarrow0,\qquad  u_{n}/\Delta_{n}^{\rho_{+}%
}\rightarrow\infty\qquad \mbox{for some }  0\leq\rho_{-}<\rho+<\tfrac{1}{2}.
\end{equation}
This condition becomes weaker when $\rho_{+}$ increases and when $\rho_{-}$ decreases.

In practice, when a Brownian motion is present, we will often translate values
of the cutoff level $u_{n}$ in terms of a number of standard deviations of the
continuous part of the semimartingale. That is, we express values of $u_{n}$
in terms of $\alpha_{n}$ where $u_{n}=\alpha_{n}(t^{-1}\int_{0}^{t}\sigma
_{s}^{2}\,ds)^{1/2}\Delta_{n}^{1/2}.$ Despite the presence of jumps, the
integrated volatility in that expression can be estimated using the small
increments of the process, since%
\begin{equation}\label{A3}
\sum_{i=1}^{[t/\Delta_{n}]}|\Delta_{i}^{n}X|^{2}1_{\{|\Delta_{i}^{n}%
X|\leq\alpha\Delta_{n}^{\varpi}\}}\stackrel{\mathbb{P}}{\longrightarrow}%
\int_{0}^{t}\sigma_{s}^{2}\,ds
\end{equation}
for any $\alpha>0$ and $\varpi\in(0,1/2).$ We can then vary the cutoff level
$\alpha_{n}$ to yield a number of (estimated)\ standard deviations of the
continuous part of the semimartingale. This data-driven choice can help
determine a range of reasonable values for the cutoff level and provide on a
path-by-path basis an equivalent, but perhaps more intuitive, scale with which
to measure the magnitude of the cutoff level $u_{n}$.

When there is no Brownian motion under the null, a different scale needs to be
used to assess the size of $u_{n}.$ For example, we can translate $u_{n}$ into
the percentage of the sample that is greater than the cutoff level, and
therefore not included in the computation of the truncated power variations.

\subsection{Testing for the presence of Brownian motion under the null}

In a first case, we set the null hypothesis to be ``the
Brownian motion is present,'' that is $\Omega_{T}^{W}$,
against the alternative $\Omega_{T}^{\mathit{noW}}\cap\Omega_{T}^{i\beta}$.

In order to construct a test, we seek a statistic with markedly different
behavior under the null and alternative. One fairly natural idea is to
consider powers less than $2$ since in the presence of Brownian motion they
would be dominated by it, while in its absence they would behave quite
differently. Specifically, the large number of small increments generated by a
continuous component would cause a power variation of order less than $2$ to
diverge to infinity. Without the Brownian motion, however, and when $p>\beta$,
the power variation converges to $0$ at exactly the same rate for the two
sampling frequencies $\Delta_{n}$ and $k\Delta_{n},$ whereas in the former
case the choice of sampling frequency will influence the magnitude of the
divergence. Taking a ratio will eliminate all unnecessary aspects of the
problem and focus on the key aspect, that of distinguishing between the
presence and absence of the Brownian motion.

Specifically, we fix a power $p\in(0,2)$ and an integer $k\geq2$, and we
consider the test statistics, which depend on $p$ and on the terminal time
$T$ and on the sequence $u_{n}$ subject to (\ref{M-7}), as follows:
\begin{equation}
S_{n} = \frac{B(p,u_{n},\Delta_{n})_{T}}{B(p,u_{n},k\Delta_{n})_{T}}.
\label{T-5}%
\end{equation}

As will become clear below, taking ratios of power variations has the
advantage of making the test statistic model-free. That is, its distribution
under the null hypothesis can be assessed without the need for the extraneous
estimation of the dynamics of the process in (\ref{1}). Obviously, these
dynamics can be quite complex with potentially jumps of various activity
levels, stochastic volatility, jumps in volatility, etc. So the fact that the
standardized test statistic can be computed without the need to estimate the
various parts of (\ref{1}) is a desirable feature. In fact, implementing the
test---that is, computing the statistic in (\ref{T-5}) and estimating its
asymptotic variance---will require nothing more than the computation of
various truncated power variations.

The first result is a law of large numbers (LLN)\ giving the probability limit
of the statistic $S_{n}$.

\begin{theorem}
\label{TT1} Under Assumption \ref{ass:A1} and if $p\in(1,2)$, we have
\begin{equation} \label{T-6}
\qquad S_{n} \stackrel{\mathbb{P}}{\longrightarrow}  \cases{
k^{1-p/2},  &\quad on the set $  \Omega_{T}^{W}$,\cr
1, & \quad  on the set $  \Omega_{T}^{\mathit{noW}}\cap\Omega_{T}^{i\beta}$,
if $ p>\beta\vee1$,  $\displaystyle\rho_{+}\leq\frac{p-1}p$.
}
\end{equation}
\end{theorem}

This result shows that, since $k^{1-p/2}>1$, for the test at hand an a
priori\thinspace reasonable critical region is $C_{n}=\{S_{n}<c_{n}\}$, for a
sequence $c_{n}$ increasing strictly to $k^{1-p/2}$: in this case the
asymptotic power is $1$ in restriction to the set described in the second
alternative above, whereas the asymptotic level depends on how fast $c_{n}$
converges to $k^{1-p/2}$.

For a more refined version of this test, with a prescribed level $a\in(0,1)$,
we need a central limit theorem (CLT)\ associated with the convergence in
(\ref{T-6}). For this we need some notation: letting $Z$ and $Z^{\prime}$ be
two independent $\mathcal{N}(0,1)$ variables, we set
\begin{equation}\label{T61}
\qquad \cases{
m_{p}=\mathbb{E}(|Z|^{p}),\cr
m_{k,p}=\mathbb{E}\bigl(|Z|^{p}\bigl|Z+\sqrt{k-1}Z^{\prime}\bigr|^{p}\bigr),\cr
\ds N(p,k)=\frac{1}{m_{2p}}  \bigl(  k^{2-p}(1+k)m_{2p}+k^{2-p}(k-1)m_{p}^{2}-2k^{3-3p/2}m_{k,p} \bigr)  .
}
\end{equation}
In terms of known functions, we have
\begin{equation}\label{T61a}
 \cases{
\ds m_{p}=\frac{2^{p/2}}{\sqrt{\pi}}\Gamma \biggl(  \frac{p+1}{2} \biggr)  ,\cr
\ds m_{k,p}=\frac{2p}{\pi}(k-1)^{p/2}\Gamma \biggl(  \frac{1+r}{2} \biggr)
^{2}F_{2,1} \biggl(  -\frac{p}{2};\frac{p+1}{2};\frac{1}{2};\frac{-1}%
{k-1} \biggr),
}
\end{equation}
where $F_{2,1}$ is Gauss's hypergeometric function (see, e.g., Section 15.1 of
\cite{abramowitzstegun}).

Then the standardized version of the CLT goes as follows (we use
$\stackrel{\mathcal{L}-(s)}{\longrightarrow}$ to denote the stable convergence
in law (see, e.g., \cite{jacodshiryaev2003} for this notion); to explain
the following statement, we recall that the convergence in law
``in restriction to a subset $\Omega_{0}$''
 is meaningless, but the stable convergence in law in restriction to
$\Omega_{0}$ makes sense):

\begin{theorem}
\label{TTCLT1} Suppose that Assumption \ref{ass:A2} holds, take $p\in(1,2)$
and let the sequence $u_{n}$ satisfy (\ref{M-7}) with $\rho_{-}>\frac
{p-1}{2p-2\beta}$. Then we have the following convergence in law:
\begin{equation}\label{T-9}
(S_{n}-k^{1-p/2})/\sqrt{V_{n}} \stackrel{\mathcal{L}-(s)}{\longrightarrow
} \mathcal{N}(0,1)\qquad   \mbox{in restriction to }\Omega_{T}^{W},
\end{equation}
where
\begin{equation}\label{T-8}
V_{n}=N(p,k) \frac{B(2p,u_{n},\Delta_{n})_{T}}{(B(p,u_{n},\Delta_{n}%
)_{T})^{2}}.
\end{equation}

\end{theorem}

We are now ready to exhibit a critical region for testing $H_{0}\dvtx\Omega
_{T}^{W}$ vs. $H_{1}\dvtx\Omega_{T}^{\mathit{noW}}\cap\Omega_{T}^{i\beta}$ using $S_{n}$
with a prescribed asymptotic level $a\in(0,1)$. Denoting by $z_{a}$ the
$a$-quantile of $N(0,1)$, that is, $\mathbb{P}(Z>z_{a})=a$ where $Z$ is
$N(0,1)$, we set
\begin{equation}\label{T-10}
C_{n} = \bigl\{S_{n}<k^{1-p/2}-z_{a}\sqrt{V_{n}}\bigr\}.
\end{equation}

\begin{theorem}\label{TT2}
Suppose that Assumption \ref{ass:A2} holds. Let $p\in(1,2)$ and
let the sequence $u_{n}$ satisfy (\ref{M-7}) with
\begin{equation} \label{ZZ3}
\frac{p-1}{2p-2\beta} = \rho_{-} < \rho_{+} = \frac{p-1}{p}
\qquad  (\mbox{hence }  p>2\beta ).
\end{equation}
Then the asymptotic level of the critical region defined by (\ref{T-10}) for
testing the null hypothesis ``the Brownian motion is
present''\ (i.e., $\Omega_{T}^{W}$ against $\Omega_{T}%
^{\mathit{noW}}\cap\Omega_{T}^{i\beta}$) equals $a$, and the asymptotic power equals
$1$.
\end{theorem}

To perform the test we need to choose $p$ and the sequence $u_{n}$. In
practice one does not know $\beta$, although it should be smaller than $1$ by
Assumption \ref{ass:A2}. Hence if we are willing to assume that $\beta$,
although unknown, is not bigger than some prescribed $\beta_{0}<1$, one should
choose $p\in(2\beta_{0},2)$, and one may take $u_{n}=\alpha  \Delta
_{n}^{\varpi}$ for some $\alpha>0$ and some $\varpi\in(0,1/2)$, and the test
can be done as soon as
\begin{equation}\label{6040}
\frac{p-1}{2p-2\beta_{0}} < \varpi <\frac{p-1}{p}.
\end{equation}
To properly separate the two hypotheses it is probably wise to choose $p$
closer to $2\beta_{0}$ than to $2$.%

\begin{figure}

\includegraphics{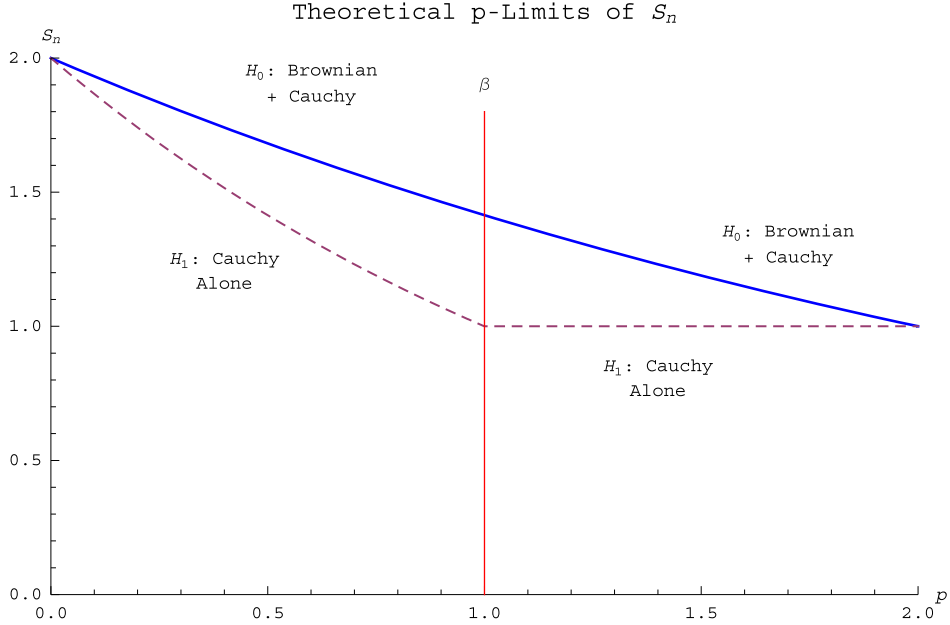}

\caption{Probability limits as a function of $p$ of the test statistic $S_{n}$
with $k=2$ in the case of a Cauchy process ($\beta=1,$ $H_{0}$) and a Brownian
plus Cauchy processes ($H_{1}$). }%
\label{fig:S4limits}%
\end{figure}

\begin{remark}\label{R1}
The first part of the consistency result (\ref{T-6}) holds also for
$p\in(0,1]$ on $\Omega_{T}^{W}$ (with basically the same proof). The second
part also holds for $\beta<p\leq1$ on the set on which $X_{t}=X_{0}%
+\sum_{s\leq t} \Delta X_{s}$ for all $t\leq T$, that is, when there is no
drift, whereas when there is a drift, $S_{n}$ converges to $k^{1-p}$ for all
$p\in(0,1]$. When $0<p\leq\beta$ the limit of $S_{n}$ is $k^{1-p/\beta}$ on
$\Omega_{T}^{\mathit{noW}}\cap\Omega_{T}^{i\beta}$when $p>1$, and also when $p\leq1$
when again there is no drift (and the proof is more involved). Figure~\ref{fig:S4limits} illustrates these various limits in the case $X$ is the sum
of a Brownian martingale plus possibly a Cauchy process (with no drift).
\end{remark}

\begin{remark}
\label{R2}The CLT necessitates $p\in(1,2)$. However, more sophisticated
techniques would allow us to prove the same result for all $p\in(0,2)$, under the
additional assumption that $\sigma_{t}$ does not vanish for $t\in\lbrack0,T]$,
on the set $\Omega_{T}^{W}$ (we still need $\beta<1$, however).
\end{remark}

\begin{remark}
\label{R4} Despite the fact that using powers less than $2$ is the most
natural way to isolate the contribution of the Brownian motion to the overall
increments of the process, it is possible to design an alternative test that
relies on powers greater than $2.$ Instead of the statistic $S_{n}$ above, we
could use the following statistic: pick $\gamma>1$ and $p^{\prime}>p>2$, and
set
\begin{equation}\label{eq:Snbar}
\overline{S}_{n} = \frac{B(p^{\prime},\gamma u_{n},\Delta_{n})_{T}%
 B(p,u_{n},\Delta_{n})_{T} B(2,\gamma u_{n},\Delta_{n})_{T}}{B(p^{\prime
},u_{n},\Delta_{n})_{T} B(p,\gamma u_{n},\Delta_{n})_{T} B(2,u_{n}%
,\Delta_{n})_{T}}. %
\end{equation}
Under Assumption \ref{ass:A1}, $\overline{S}_{n}$ converges in probability to
$\gamma^{p^{\prime}-p}$ on the set $\Omega_{T}^{W}\cap\Omega_{T}^{i\beta}$,
and to $\gamma^{p^{\prime}-p+2-\beta}$ on the set $\Omega_{T}^{\mathit{noW}}\cap
\Omega_{T}^{i\beta}$, as soon as $\rho_{+}<\frac{p-2}{2p-2\beta}$. We also
have a CLT under Assumption \ref{ass:A2} and if $\rho_{+}<\frac{2p-4}{11p-10}%
$. Under $H_{0},$ $\overline{S}_{n}$ is model-free, just like $S_{n}$ is. So
one can, in an obvious way, construct a test based on $\overline{S_{n}}$ and
which satisfies the claims of Theorem \ref{TT2}, under suitable conditions on
the cutoff levels $u_{n}$. However, simulations studies suggest that the
statistic $\overline{S_{n}}$ is not as well behaved as $S_{n}$, and so we do
not pursue its study further.
\end{remark}

\subsection{Testing for the absence of Brownian motion under the null}

In a second case, we set the null hypothesis to be ``the
Brownian motion is absent,'' that is, $\Omega_{T}^{\mathit{noW}}%
\cap\Omega_{T}^{i\beta}$. Designing a test under this null is trickier because
the model becomes a pure jump (plus perhaps a drift) process, and we are
aiming for a test that remains model-free even for this model. That is, we are
looking for a statistic whose limiting behavior under the null, despite being
driven by what is now a pure jump process, does not depend on the
characteristics of the pure jump process, such as its degree of activity
$\beta$, since those characteristics are a priori unknown.

This can be achieved as follows. We choose a real $\gamma>1$ and a sequence
$u_{n}$ satisfying (\ref{M-7}) and define the test statistic
\begin{equation}\label{T-11}
S_{n}^{\prime} = \frac{B(2,\gamma u_{n},\Delta_{n})_{T} U(u_{n},\Delta
_{n})_{T}}{B(2,u_{n},\Delta_{n})_{T} U(\gamma u_{n},\Delta_{n})_{T}}.
\end{equation}
To understand the construction of this test statistic, recall that in a power
variation of order $2$ the contributions from the Brownian and jump components
are of the same order. But once the power variation is properly truncated, the
Brownian motion will dominate it if it is present. And the truncation can be
chosen to be sufficiently loose that it retains essentially all the increments
of the Brownian motion at cutoff level $u_{n}$ and a fortiori $\gamma u_{n}$,
thereby making the ratio of the two truncated quadratic variations converge to
$1$ under the alternative hypothesis. On the other hand, if the Brownian
motion is not present, then the nature of the tail of jump distributions is
such that the difference in cutoff levels between $u_{n}$ and $\gamma u_{n}$
remains material no matter how far we go in the tail, and the limit of that
same ratio will reflect it: it will now be $\gamma^{2-\beta}$ under
assumptions made specific in the formal theorems below. Since absence of a
Brownian motion is now the null hypothesis, the issue is then that this limit
depends on the unknown $\beta.$

Canceling out that dependence is the role devoted to the ratio of the number
of large increments, the $U$'s, in (\ref{T-11}). The $U$'s are
always dominated by the jump components of the model whether the Brownian
motion is present or not. Their inclusion in the statistic is merely to ensure
that the statistic is model-free, by effectively canceling out the dependence
on the jump characteristics that emerges from the ratio of the truncated
quadratic variations. Indeed, the limit of the ratio of the $U$'s is
$\gamma^{\beta}$ under both the null and alternative hypotheses. As a result,
the probability limit of $S_{n}^{\prime}$ will be $\gamma^{2}$ under the null,
independent of $\beta$.

Our first result states this precisely, establishing the limiting behavior of
the statistic in terms of convergence in probability:

\begin{theorem}
\label{TT3} Let the sequence $u_{n}$ satisfy (\ref{M-7}), and suppose that
Assumption~\ref{ass:A1} holds. Then
\begin{equation}
\label{T-12}S^{\prime}_{n} \stackrel{\mathbb{P}}{\longrightarrow}  \cases{
\gamma^{2},  & \quad\mbox{on the set }   $\Omega_{T}^{\mathit{noW}}\cap\Omega
_{T}^{i\beta}$,\cr
\gamma^{\beta}, & \quad\mbox{on the set }   $\Omega_{T}^{W}\cap\Omega_{T}^{i\beta
}$.
}
\end{equation}

\end{theorem}

For a test with a prescribed level we need a standardized CLT.

\begin{theorem}
\label{TTCLT2} Suppose that Assumption \ref{ass:A1} holds with $\beta
^{\prime\prime}<\frac{\beta}{2+\beta}$ and $\beta^{\prime}<\frac{\beta}{2}$,
and (\ref{M-7}) holds with $\rho_{+}<\frac{1}{2+\beta}\wedge\frac{2}{5\beta
}\wedge\frac{2-\beta}{3\beta}$. Then we have
\begin{equation}\label{eq:CLT_S'}
(S_{n}^{\prime}-\gamma^{2})/\sqrt{V_{n}^{\prime}} \stackrel{\mathcal{L}%
-(s)}{\longrightarrow} \mathcal{N}(0,1)\qquad  \mbox{in restriction to
}\Omega_{T}^{\mathit{noW}}\cap\Omega_{T}^{i\beta},
\end{equation}
where $V_{n}^{\prime}$ is given by the following formula:
\begin{eqnarray} \label{CLT V'}
V_{n}^{\prime}  &  =&\gamma^{4} \biggl(  \frac{B(4,u_{n},\Delta_{n})_{T}%
}{(B(2,u_{n},\Delta_{n})_{T})^{2}}+\frac{1}{U(u_{n},\Delta_{n})}
\nonumber\\[-8pt]\\[-8pt]
&&\hphantom{\gamma^{4} \biggl(}{}    + \biggl(  1-\frac{2}{\gamma^{2}} \biggr)   \biggl(  \frac
{B(4,\gamma u_{n},\Delta_{n})_{T}}{(B(2,\gamma u_{n},\Delta_{n})_{T})^{2}%
}+\frac{1}{U(\gamma u_{n},\Delta_{n})} \biggr)   \biggr) .\nonumber
\end{eqnarray}
\end{theorem}

Hence a critical region for testing $H_{0}\dvtx\Omega_{T}^{\mathit{noW}}\cap\Omega
_{T}^{i\beta}$ vs. $H_{1}\dvtx\Omega_{T}^{W}\cap\Omega_{T}^{i\beta}$ is
\begin{equation}
\label{T-113}C^{\prime}_{n} = \bigl\{S_{n}^{\prime}<\gamma^{2}-z_{a}\sqrt
{V_{n}^{\prime}}\bigr\}.
\end{equation}

\begin{theorem}
\label{TT4} Suppose that Assumption \ref{ass:A1} holds with $\beta
^{\prime\prime}<\frac{\beta}{2+\beta}$ and $\beta^{\prime}<\frac{\beta}2$, and
(\ref{M-7}) holds with $\rho_{+}<\frac1{2+\beta}\wedge\frac2{5\beta}
\wedge\frac{2-\beta}{3\beta}$. Then the asymptotic level of the critical
region $C^{\prime}_{n}$ defined by (\ref{T-113}) for testing the null
hypothesis ``the Brownian motion is absent''
 (i.e., $\Omega_{T}^{\mathit{noW}}\cap\Omega_{T}^{i\beta}$ against $\Omega_{T}%
^{W}\cap\Omega_{T}^{i\beta}$) equals $a$, and the asymptotic power equals $1$.
\end{theorem}

If we take again $u_{n}=\alpha \Delta_{n}^{\varpi}$, the test can be
performed if $\alpha>0$ and
\begin{equation}\label{603}
0 <\varpi < \frac{1}{2+\beta}\wedge\frac{2}{5\beta}\wedge\frac{2-\beta}%
{3\beta}
\end{equation}
(always smaller than $1/2$). This requirement is constraining, because $\beta$
is unknown, and may typically be close to $2$ if we believe in the null
hypothesis. Therefore in practice we must assume that $\beta$ does not exceed
a given $\beta_{0}\in\lbrack1,2)$. This means that this limiting index
$\beta_{0}$ is given a priori, and we do the test under the Assumption
\ref{ass:A1} with $2\beta^{\prime}<\beta\leq\beta_{0}$ and $\beta
^{\prime\prime}<\frac{\beta}{2+\beta}$, with $\varpi$ subject to the
(feasible) condition
\begin{equation}\label{604}
0 <\varpi <\frac{2-\beta_{0}}{3\beta_{0}}. %
\end{equation}
These facts are not really surprising: first, by (\ref{T-12}) we know that the
statistic $S_{n}^{\prime}$ properly separates the two hypotheses only when
$\beta$ is not too close to $2$. And, second, when $\beta$ becomes very close
to $2$, the paths of $X$ have big jumps but also the compensated sum of
small\ jumps looks more and more like a Brownian path, even on the set
$\Omega_{T}^{\mathit{noW}}$.

\begin{remark}
\label{R5} It is possible to design an alternative statistic with similar
properties but make no use of the $U$'s. Instead of the statistic
$S_{n}^{\prime}$ in (\ref{T-11}), we could use the following statistic: pick
$\gamma>1,$ $\kappa\geq1$ and $p>2$, and set
\begin{equation}\label{eq:Snbarprime}
\overline{S}{}_{n}^{\prime} = \frac{B(2,u_{n},\Delta_{n})_{T}
B(p,\kappa\gamma u_{n},\Delta_{n})_{T} }{B(2,\gamma u_{n},\Delta_{n}%
)_{T} B(p,\kappa u_{n},\Delta_{n})_{T}}.
\end{equation}
Under Assumption \ref{ass:A1}, $\overline{S}{}_{n}^{\prime}$ converges
in probability to $\gamma^{p-2}$ on the set
$\Omega_{T}^{\mathit{noW}}\cap\Omega _{T}^{i\beta}$, and to
$\gamma^{p-\beta}$ on the set $\Omega_{T}^{W}\cap \Omega_{T}^{i\beta}$,
as soon as $\rho_{+}<\frac{p-2}{2p}$. The ratio of $p$th power
variations plays a similar role to that of the $U$'s, namely to cancel
out the dependence of the $p$-lim of $\overline{S}{}_{n}^{\prime}$ on
$\beta$ under the null. The fixed scaling factor $\kappa$ allows us to
use different cutoff levels for the two powers $p$ and $2$ without
affecting the probability limit of the statistic. We also have a CLT if
$\rho_{+}<\frac{2-\beta}{3\beta}$. Under $H_{0},$
$\overline{S}{}_{n}^{\prime}$ is model-free, just like $S_{n}^{\prime}$
is, and so a test follows. But as was the case for the statistic
$\overline{S}_{n}$ proposed in (\ref{eq:Snbar}), simulations studies
suggest that $\overline{S}{}_{n}^{\prime}$ is not as well behaved as
$S_{n}$.
\end{remark}

\begin{remark}
\label{R10} In Theorems \ref{TTCLT1} and \ref{TT2} the rate of convergence is
hidden because of the standardization, but it is $1/\sqrt{\Delta_{n}}$,
clearly optimal since there are $1+[T/\Delta_{n}]$ observation altogether. In
Theorems \ref{TTCLT2} and \ref{TT4} the rate is $1/u_{n}^{\beta/2}$, which is
again ``optimal'' when we only use the
increments bigger than $u_{n}$ [more precisely, if we were able to observe
exactly all jumps of $X$ with size bigger than $u_{n}$, this rate would be the
optimal one, up to a $\log(1/u_{n})$ term]. However, for those theorems we also
have to choose $u_{n}$: the smallest $u_{n}$ is, compared to $\Delta_{n}$, the
biggest the actual rate is, but we are limited in this choice by the upper
bound on $\rho_{+}$. For example if we take $u_{n}=\alpha\Delta_{n}^{\varpi}$,
and due to (\ref{604}), the best rate is ``almost''\ $1/\Delta_{n}^{\beta(2-\beta_{0})/6\beta_{0}}$.
\end{remark}

\section{Simulation results}\label{sec:MC}

We now report simulation results documenting the finite sample performance of
the test statistics $S_{n}$ and $S_{n}^{\prime}$. We calibrate the values to
be realistic for a liquid stock trading on the NYSE, and we consider an
observation length of $T=21$ days (one month)\ sampled every five seconds.%

We conduct simulations to determine the small sample behavior of the two
statistics $S_{n}$ and $S_{n}^{\prime}$ under their respective null and
alternative hypotheses. The tables and graphs that follow report the results
of $5000$ simulations. The data generating process is the stochastic
volatility model $dX_{t}=\sigma_{t}\,dW_{t}+\theta\,dY_{t},$ with $\sigma
_{t}=v_{t}^{1/2}$, $dv_{t}=\xi(\eta-v_{t})\,dt+\phi v_{t}^{1/2}\,dB_{t}+dJ_{t}$,
$\mathbb{E}[dW_{t}\,dB_{t}]=\rho\, dt$, $\eta^{1/2}=0.25,$ $\phi=0.5,$ $\xi=5$,
$\rho=-0.5,$ $J$ is a compound Poisson jump process with jumps that are
uniformly distributed on $[-30\%,30\%]$ and $X_{0}=1$. The jump process $Y$ is
a $\beta$-stable process with $\beta=1$, that is, a Cauchy process (which has
infinite activity, and will be our model under $\Omega_{T}^{i\beta}$; this is
a borderline case for the statistics $S_{n}$ under the null, nevertheless we
will see that this statistic behaves well). Given $\eta,$ the scale parameter
$\theta$ (or equivalently $A$) of the stable process in simulations is
calibrated to deliver different various values of the tail probability
$P=\mathbb{P}(|\Delta Y_{t}|\geq4\eta^{1/2}\Delta_{n}^{1/2})$. In the various
simulations' design, we hold $\eta$ fixed. Therefore the tail probability
parameter $P\ $controls the relative scale of the jump component of the
semimartingale relative to its continuous counterpart. We set $\theta$ such
that neither of the two components of the model, $\sigma_{t}\,dW_{t}$ and
$\theta Y_{t},$ is negligible compared to the other when the hypothesis states
that they should  both be present. We achieve this by computing the expected
percentage of the total quadratic variation attributable to jumps on a given
path from the model, and set it to values that range from $5\%$ and
$95\%.$\looseness=-1

\begin{figure}[b]

\includegraphics{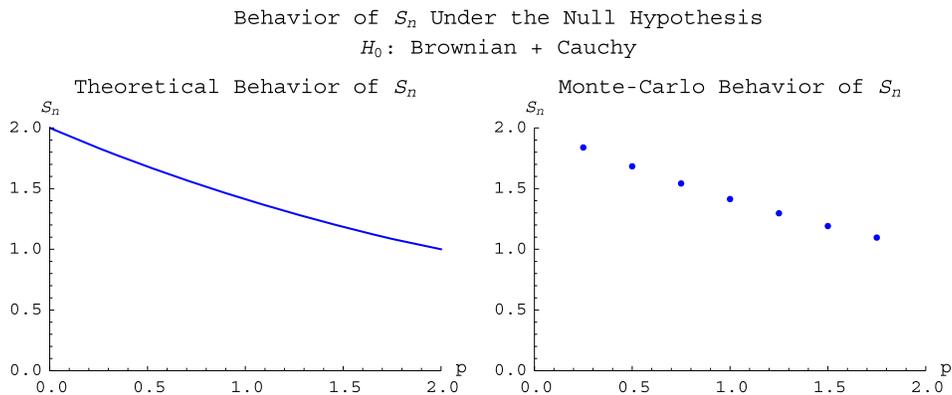}

\caption{Theoretical and Monte Carlo behavior of $S_{n}$ as a function of the
power $p$ under the null hypothesis where a Brownian motion is present, in
addition to a pure jump (Cauchy) process.}%
\label{fig:S4H0Curves}%
\end{figure}

\subsection{The first test}

\begin{figure}

\includegraphics{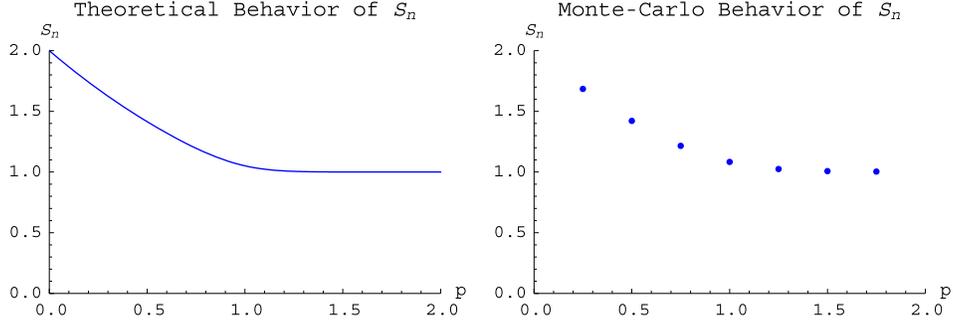}

\caption{Theoretical and Monte Carlo behavior of $S_{n}$ as a function of the
power $p$ under the alternative hypothesis where a Brownian is absent.}\label{fig:S4H1Curves}%
\end{figure}

The statistic $S_{n}$ is implemented with $k=2$ and values of $p$ that range
from $0$ to $2$ (recall Remark \ref{R2}). Figure \ref{fig:S4H0Curves} compares
the theoretical and Monte Carlo behavior of $S_{n}$ as a function of the power
$p$ under the null hypothesis where a Brownian motion is present, in addition
to a Cauchy pure jump process. Figure \ref{fig:S4H1Curves} shows the
corresponding results under the alternative hypothesis, where there is no
Brownian motion. The theoretical curves are computed from the expected values
of the truncated power variations using the exact density of the increments at
the sampling interval $\Delta_{n}=5$ seconds, rather than their asymptotic
limits for $\Delta_{n}\rightarrow0.$ This introduces a slight Jensen's
inequality effect in the figure but appears to capture well the small sample
behavior of the statistic.

Recall that for concreteness $\alpha$ is expressed as a number of standard
deviations of the Brownian part of $X$: that is, the level of truncation
$u_{n}$ is expressed in terms of the number $\alpha$ of standard deviations of
the continuous martingale part of the process, defined in multiples of the
long-term volatility parameter $\eta^{1/2}$: $ \alpha$ is defined by
$u_{n}=\alpha\eta^{1/2}\Delta_{n}^{1/2}.$ Our view of the joint choice of
$(\varpi,\alpha)$ is that they are not independent parameters in finite
sample: they are different parameters for asymptotic purposes but in finite
samples the only relevant quantity is the actual resulting cutoff size $u_{n}%
$. This is why we are reporting the values of the cutoffs $u_{n}$ in the form
of the $\alpha$ that would correspond to $\varpi=1/2.$ This has the advantage
of providing an easily interpretable size of the cutoff compared to the size
of the increments that would be expected from the Brownian component of the
process: we can then think in terms of truncating at a level that corresponds
to $\alpha=4$, $6,$ etc., standard deviations of the continuous part of the
model. Since the ultimate purpose of the truncation is either to eliminate or
conserve that part, it provides an immediate and intuitively clear reference
point. Of course, given $u_{n}$ and this $\alpha,$ it is possible to back this
into the value of the $\alpha$ corresponding to any $\varpi,$ for that given
sample size, including the value(s) of $\varpi$ that satisfy the required
inequalities imposed by the asymptotic results. This approach would lose its
effectiveness if we were primarily interested in testing the validity of the
asymptotic approximation as the sample size varies, but for applications, by
definition on a finite sample, it seems to us that the interpretative
advantage outweighs this disadvantage.

\begin{table}[b]
\caption{Testing $ H_{0}\dvtx\Omega_{T}^{W}$ vs. $H_{1}\dvtx\Omega
_{T}^{\mathit{noW}}\cap\Omega_{T}^{i\beta}$:
Monte Carlo rejection rate for the test for the presence of a Brownian motion
using the statistic $S_n$}\label{table:MCS4H0}
\begin{tabular*}{\textwidth}{@{\extracolsep{\fill}}ld{3.0}cccd{2.1}cd{2.1}c@{}}
\hline
&& \multicolumn{7}{c@{}}{\textbf{Sample rejection rate (\%) for power} $\bolds{p}$}\\[-5pt]
&&\multicolumn{7}{c@{}}{\hrulefill}\\
\multirow{2}{50pt}[12pt]{\textbf{Degree of truncation}
$\bolds{\alpha}$} & \multicolumn{1}{c}{\multirow{2}{59pt}[12pt]{\centering\textbf{Test theoretical level}}}     &
\textbf{0.25} & \textbf{0.5} & \textbf{0.75} & \textbf{1.0} & \textbf{1.25} &
\textbf{1.5} & \textbf{1.75}\\
\hline
6 & 10\%   & 9.1 & 9.4 & 9.4 & 9.2 & 9.1 & 9.3 & 8.9\\
& 5\% &   4.6 & 4.7 & 4.8 & 4.7 & 4.5 & 4.3 & 4.1\\
7 & 10\%   & 9.7 & 9.7 & 9.8 & 9.8 & 9.9 & 10.2 & 9.8\\
& 5\% &   5.0 & 5.0 & 5.1 & 5.0 & 4.9 & 4.7 & 4.4\\
8 & 10\%   & 9.7 & 9.9 & 9.9 & 10.0 & 9.9 & 10.1 & 9.9\\
& 5\% &   5.0 & 5.1 & 5.1 & 5.1 & 4.9 & 4.8 & 4.5\\
\hline
\end{tabular*}%
\end{table}

The statistic in the plots is computed with a truncation level corresponding
to $\alpha=7$. Table \ref{table:MCS4H0} looks at the dependence of the results
on the choice of $\alpha.$

\begin{figure}

\includegraphics{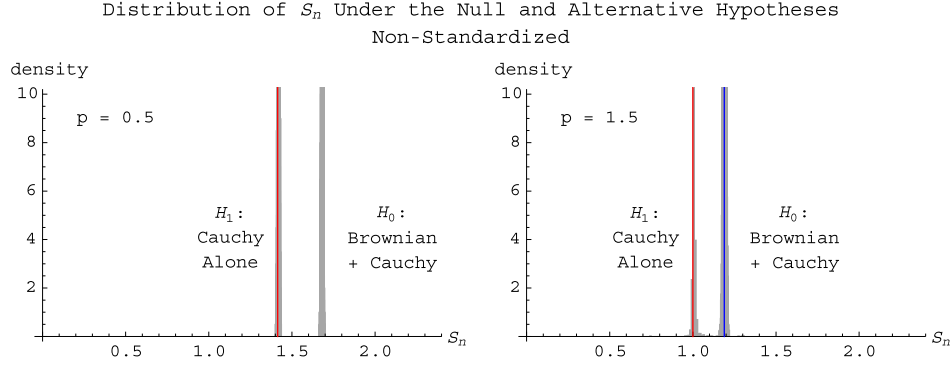}

\caption{Nonstandardized distribution of $S_{n}$ under the null and
alternative hypotheses for two values of $p.$}%
\label{fig:S4H0H1HistNonStd}%
\end{figure}

Next, we report in Figure \ref{fig:S4H0H1HistNonStd} histograms of the values
of the unstandardized $S_{n}$ computed under $H_{0}\dvtx\Omega_{T}^{W}$ and
$H_{1}\dvtx\Omega_{T}^{\mathit{noW}}\cap\Omega_{T}^{i\beta}$, respectively, and with the
same level of truncation $\alpha=7$. The vertical lines represent the
anticipated limits of the statistic in the two situations, $k^{1-p/2}$ under
$H_{0}$ and either $1$ when $p>\beta$ or $k^{1-p/\beta}$ when $0<p\leq\beta$
under $H_{1}$, based on Theorem \ref{TT1} and Remark \ref{R1}. Since here
$\beta=1$, the two graphs with $p=0.5$ and $p=1.5$ illustrate the two
situations where $p<\beta$ and $p\geq\beta.$

\begin{figure}[b]

\includegraphics{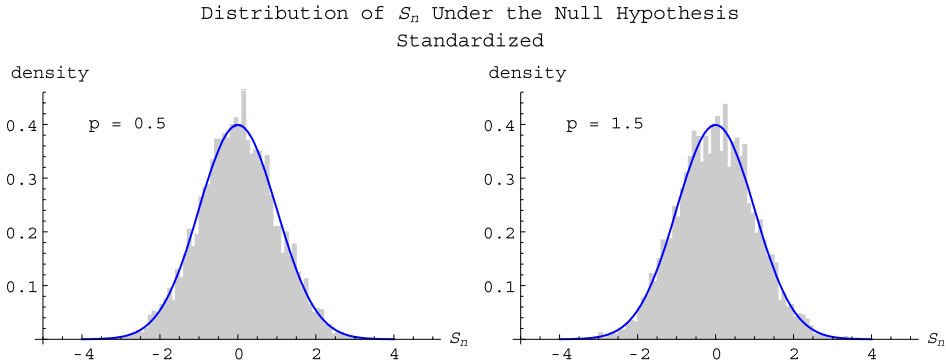}

\caption{Standardized distribution of $S_{n}$ under the null hypothesis of a
Brownian motion present for two values of $p.$ The histogram represents the
small sample distribution while the solid curve is the asymptotic
$\mathcal{N}(0,1)$ density.}
\label{fig:S4H0HistStd}
\end{figure}

Figure \ref{fig:S4H0HistStd} reports the Monte Carlo distribution of the
statistic $S_{n},$ standardized according to Theorem \ref{TTCLT1}, compared to
the limiting $\mathcal{N}(0,1)$\break distribution. Table \ref{table:MCS4H0}
reports the Monte Carlo rejection rates of the test of $H_{0}\dvtx\Omega_{T}^{W}$
vs. $H_{1}\dvtx\Omega_{T}^{\mathit{noW}}\cap\Omega_{T}^{i\beta}$ at the $10\%$ and $5\%$
level, using the test statistic $S_{n},$ for various levels of truncation
$\alpha$. We find that the test behaves well, with empirical test levels close
to their theoretical counterparts.

\subsection{The second test}

We now turn to the second problem, that of testing $H_{0}\dvtx\Omega_{T}^{\mathit{noW}}%
\cap\Omega_{T}^{i\beta}$ vs. $H_{1}\dvtx\Omega_{T}^{W}\cap\Omega_{T}^{i\beta}$.
For this test, $S_{n}^{\prime}$ is implemented with a second truncation level
twice as large as the first, that is, $\gamma=2.$ The simulation evidence
suggests that the results are largely similar for values of $\gamma$ within a
range of $1.5$ to $2.5$.\ Parameter values are identical to those employed for
the first test. Since there is no Brownian motion under the null, the
truncation level $u_{n}$ is set in terms of the percentage of observations
that are excluded by the truncation. For comparison with the truncation levels
employed in the first test, we report it here again in terms of $\alpha,$ a
number of standard deviations for the Brownian motion using the same parameter
values as under the first test's null, or this test's alternative
hypothesis.\looseness=1

\begin{figure}[b]

\includegraphics{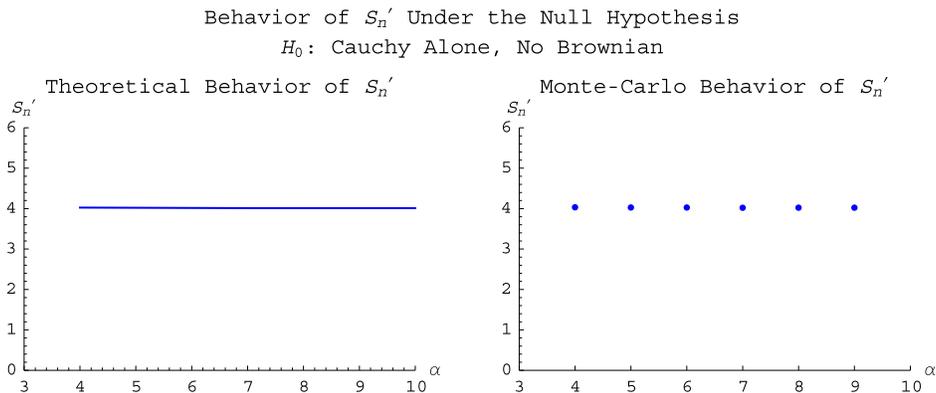}

\caption{Theoretical and Monte Carlo behavior of $S_{n}^{\prime}$ as a
function of the truncation level $\alpha$ under the null hypothesis where a
Brownian motion is absent. The model is a pure jump (Cauchy) process.}%
\label{fig:S5H0Curves}%
\end{figure}

\begin{figure}

\includegraphics{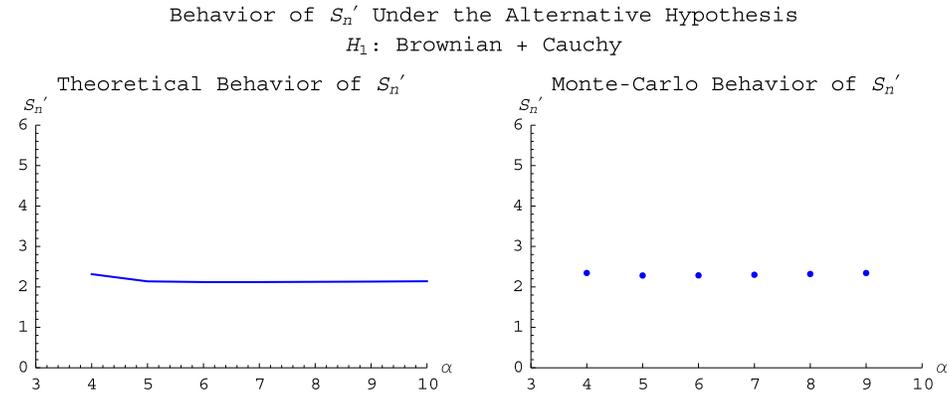}

\caption{Theoretical and Monte Carlo behavior of $S_{n}^{\prime}$ as a
function of the truncation level $\alpha$ under the alternative hypothesis
where a Brownian motion is present, in addition to a pure jump
(Cauchy) process.}\label{fig:S5H1Curves}
\end{figure}

Under the null, the model is driven exclusively by the Cauchy process.
Figure~\ref{fig:S5H0Curves} shows the limiting value of $S_{n}^{\prime}$ under
$H_{0},$ as a function of the truncation level $\alpha,$ comparing the
theoretical limit of $\gamma^{2}=4$ given in Theorem \ref{TT3} (left graph)
and the corresponding average value of $S_{n}^{\prime}$ from the Monte Carlo
simulations (right graph). Figure \ref{fig:S5H1Curves} shows the corresponding
values under the alternative hypothesis, where the increments of $X$ are now
generated by a Brownian motion plus a Cauchy process. The theoretical limit on
the left graph is computed from the expected values under the exact
distribution of the increments at the sampling frequency $\Delta_{n}$ rather
than the $p$-lim $\gamma^{2-\beta}=2$ obtained in the limit where $\Delta
_{n}\rightarrow0,$ with the same remark about Jensen's inequality applying
here. We note that for small truncation levels ($\alpha=4$)\ the interaction
of the Brownian and the stable processes is material, driving the actual limit
above $2.$ If desired, small sample corrections for this interaction can be
implemented along the same lines as in Section 5 of \cite{yacjacod09b}.

The test statistic in simulations under the alternative appears to be slightly
biased upwards. Quite naturally, this effect worsens as the pure jump process
gets closer to a Brownian motion (for instance if $\beta=1.5$ instead of $0.5$
or $1$), and/or when the scale parameter $\theta$ of the jump process
increases since that makes isolating the effect of the Brownian motion
component of the model relatively more difficult.

Generally speaking, $S_{n}^{\prime}$ is, under its alternative, more finicky
than $S_{n}$ is under either its null or alternative. The reason for this is
that $S_{n}^{\prime}$ requires under $H_{1}$ a Goldilocks-like conjunction of
factors whereby the Brownian motion component of the model is sufficiently
large to drive the behavior of the ratio of truncated quadratic variations,
while the jump component of the model cannot be so small as to render
inaccurate the ratio of the number of increments larger than the truncation level.

\begin{figure}

\includegraphics{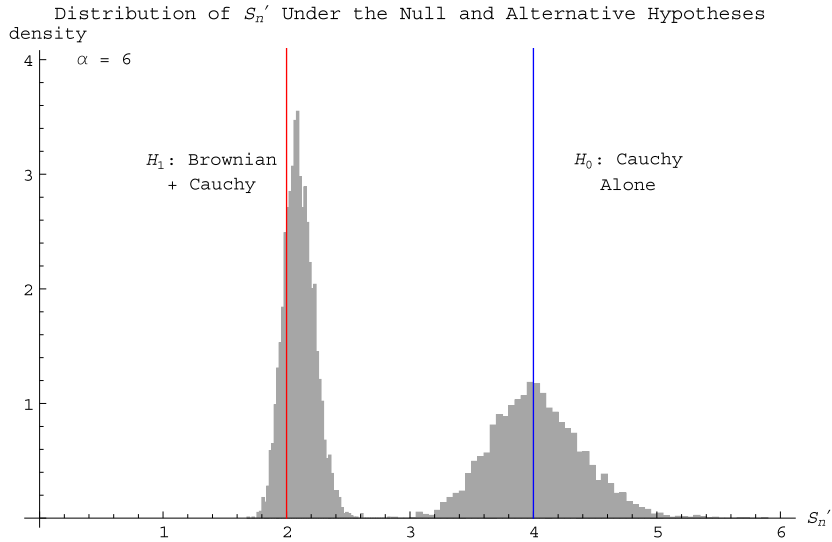}

\caption{Nonstandardized distribution of $S_{n}^{\prime}$ under the null and
alternative hypotheses.}\label{fig:S5H0H1HistNonStd}
\end{figure}

\begin{figure}[b]

\includegraphics{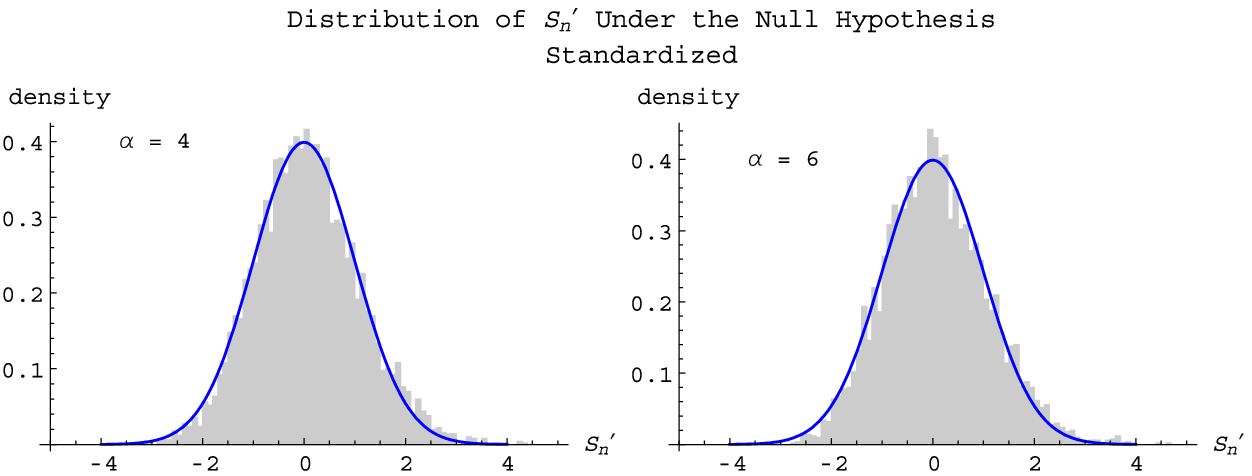}

\caption{Standardized distribution of $S_{n}^{\prime}$ under the null
hypothesis of a Brownian motion absent for two values of the truncation level.
The histogram represents the small sample distribution while the solid curve
is the asymptotic $\mathcal{N}(0,1)$ density.}\label{fig:S5H0HistStd}
\end{figure}

Figure \ref{fig:S5H0H1HistNonStd} reports the Monte Carlo distributions of
$S_{n}^{\prime}$ under $H_{0}$ and $H_{1};$ the vertical lines represent the
theoretical limits. Under $H_{1},$ we note again that $S_{n}^{\prime}$ is
slightly biased upwards. Fortunately, this bias is limited to $H_{1}$ so it
does not adversely affect the implementation of the test per se, which is
based on the behavior of $S_{n}^{\prime}$ under $H_{0}$. But it can affect the
interpretation of the results of the test implemented on real data, since, as
we will see below, we will find empirical values of $S_{n}^{\prime}$
below~$4$. Figure \ref{fig:S5H0HistStd} reports the Monte Carlo and asymptotic
distribution of the statistic $S_{n}^{\prime}$ standardized under $H_{0}$ as
prescribed by Theorem \ref{TTCLT2}.

As said above, the histograms are computed using $T=21$ days (one
month)\ sampled every five seconds. With this length of the series, the
empirical distribution of the statistic is very well approximated by its
asymptotic $\mathcal{N}(0,1)$ limit. Shorter time periods (such as $T=1$ day)
tend to result in right-skewness of the Monte Carlo distribution of
$S_{n}^{\prime}$. We do not view the need for a longer series as a serious
obstacle to the empirical implementation of the test since one would not
typically expect the Brownian motion component of the model to be turned on or
off on a daily basis: one would expect the market to operate in such a way
that the Brownian component is either there all the time or not there at all.
But if an answer is nevertheless desired on a day-by-day basis, then the first
test can always be implemented, as it requires substantially shorter time spans.

\begin{table}
\tablewidth=260pt
\caption{Testing $
H_{1}\dvtx\Omega_{T}^{\mathit{noW}}\cap\Omega_{T}^{i\beta}$
vs. $ H_{0}\dvtx\Omega_{T}^{W}$:
Monte Carlo rejection rate for the test for the absence of a Brownian
motion using the statistic $S_{n}%
^{\prime}$}\label{table:MCS5H0}
\begin{tabular*}{260pt}{@{\extracolsep{4in minus 4in}}lccccc@{}}
\hline
  & \multicolumn{5}{c@{}}{\textbf{Sample rejection rate (\%)}}\\
 & \multicolumn{5}{c@{}}{\textbf{for truncation level} $\bolds{\alpha}$}\\[-5pt]
&\multicolumn{5}{c@{}}{\hrulefill}\\
\textbf{Test theoretical level}&   $\bolds{4}$ & $\bolds{5}$ & $\bolds{6}$ & $\bolds{7}$ & $\bolds{8}$\\
\hline
10\%  & $9.0$ & $9.2$ & $9.1$ & $9.0$ & $9.2$\\
\hphantom{0}5\%   & $4.0$ & $4.2$ & $4.1$ & $4.2$ & $4.2$\\
\hline
\end{tabular*}
\end{table}

Finally, the test's rejection rate under the null hypothesis is reported in
Table \ref{table:MCS5H0}. Since the test is one-sided (we reject $H_{0}$ when
the standardized $S_{n}^{\prime}$ is too low), the right-skewness of the
statistic visible in Figure \ref{fig:S5H0HistStd} results in a slight
under-rejection by the test.

\section{Empirical results}\label{sec:data}

In this section, we apply the two test statistics to real data, consisting of
all transactions recorded during the year 2006 on two of the most actively
traded stocks, Intel (INTC) and Microsoft (MSFT). The data source is the TAQ
database. Using the correction variables in the dataset, we retain only
transactions that are labeled ``good trades''
 by the exchanges: regular trades that were not corrected, changed, or
signified as cancelled or in error; and original trades which were later
corrected, in which case the trade record contains the corrected data for the
trade. Beyond that, no further adjustment to the raw data are made.

\begin{figure}

\includegraphics{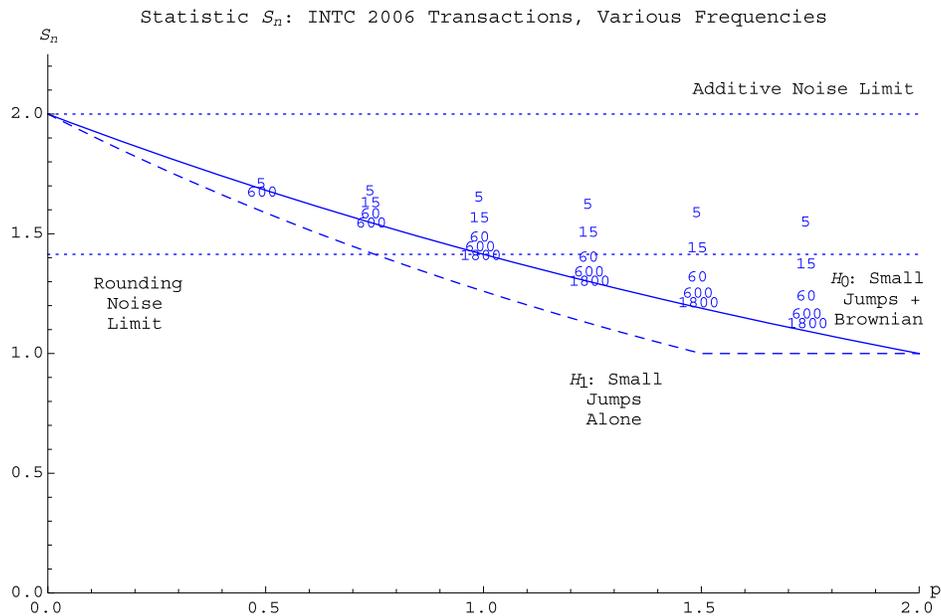}

\caption{Empirical estimates of $S_{n}$ at various values of $p$ and sampling
frequencies from all Intel transactions during 2006.}\label{fig:S4INTC}
\end{figure}

\begin{figure}

\includegraphics{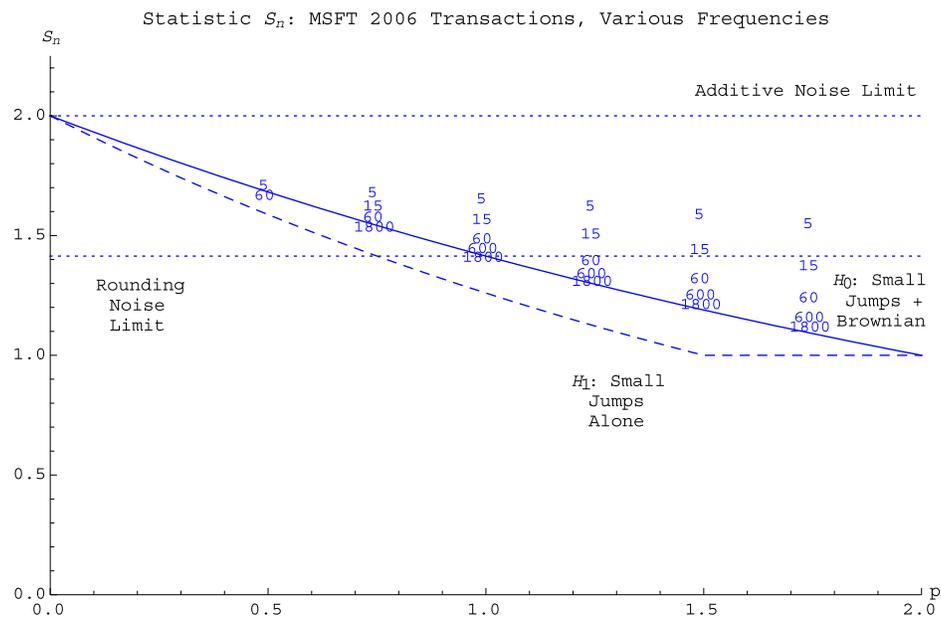}

\caption{Empirical estimates of $S_{n}$ at various values of $p$ and sampling
frequencies from all Microsoft transactions during 2006.}\label{fig:S4MSFT}
\end{figure}

We first consider the test where the null hypothesis consists of a continuous
component being present. Figures \ref{fig:S4INTC} and \ref{fig:S4MSFT} show
the values of the test statistic $S_{n},$ plotted for a range of values of the
power $p,$ for the two data series. The empirical values of $S_{n}$ are
labeled on the plots with numbers representing the sampling interval employed,
in seconds, with values ranging from $\Delta_{n}=5$ seconds to $\Delta_{n}=30$
minutes. In addition to the empirical estimates, the figures display the two
limits of $S_{n}$ under the null where a Brownian is present and the
alternative hypothesis where it is absent. The theoretical limits correspond
to those given in Figures \ref{fig:S4H0Curves} and \ref{fig:S4H1Curves},
except that the theoretical limit under $H_{1}$ (no Brownian present)\ is
plotted for a value of $\beta=1.6$, in line with the estimates of $\beta$
given in \cite{yacjacod09b} for these data series. Quite naturally, the closer
$\beta$ is to $2,$ the closer the jump component can mimic the behavior of a
Brownian motion and the harder it becomes to tell the two hypotheses apart.
The limit under $H_{0}$ is independent of $\beta.$ Also on the figures are the
two limits corresponding to the situation where market microstructure noise
dominates. We include the two polar cases where the noise is either of a pure
additive form or of a pure rounding form.

When the observations are blurred with either an additive white noise or with
noise due to rounding,\ the respective limits are then $2$ and $\sqrt{2}$.
Indeed, suppose that instead of observing the exact value of $X$ we have on
top of it an additive white noise, that is we observe $X_{i\Delta_{n}}+Z_{i}$
(at stage $n$, where the $Z_{i}$'s are i.i.d., and independent of the process
$X$). If we suppose that $Z_{i}$ has a density which is continuous and nonvanishing at $0$, then the noise is the leading factor in the behavior of
$B(p,u_{n},\Delta_{n})_{T}$ as soon as $p/(2(p+1))\geq\rho_{+}$ [recall
(\ref{M-7})]. In this case, the variables $(\Delta_{n}/u_{n}^{p+1}%
) B(p,u_{n},\Delta_{n})_{T}$ converge in probability to $TC_{p}$ for some
constant $C_{p}$, and thus $S_{n}$ converges in probability to the sampling
frequency ratio $k,$ which is $2$ here. When the noise is pure rounding\ at
some level $\alpha_{n}$, then again it is the leading factor and $\sqrt
{\Delta_{n}} \alpha_{n}^{1-p} B(p,u_{n},\Delta_{n})_{T}$ converges in
probability to some positive limiting variable, as soon as $\alpha_{n}%
/u_{n}\rightarrow0$ and $\alpha_{n}^{2}\Delta_{n}\rightarrow\infty$. Thus
$S_{n}$ converges in probability to $\sqrt{k}$ [when $\alpha_{n}>u_{n}$ we
have $B(p,u_{n},\Delta_{n})_{T}=0$ and then $S_{n}$ is not even well defined;
however, here the truncation level $u_{n}$ used in practice is quite bigger
than the rounding level of $1$ cent].

The values of $\alpha$ are similar to those employed in simulations, and
indexed in terms of standard deviations of the continuous martingale part of
the log-price: we first estimate the volatility of the continuous part of $X$
using the small increments, those of order $\Delta_{n}^{1/2}$, and then use
that estimate to form the cutoff level used in the construction of the test
statistic. To account for potential time series variation in the volatility
process $\sigma_{t}$, that procedure is implemented separately for each day
and we compute the sum, for that day, of the absolute value of the increments
that are smaller than the cutoff, to the appropriate power $p$. For the full
year, we then add the truncated power variations computed for each day.

The results in both Figures \ref{fig:S4INTC} and \ref{fig:S4MSFT} tell a
similar story. First, the empirical estimates are always on the side away from
the limit under $H_{1},$ indicating that the null hypothesis of a Brownian
motion present will not be rejected. Second, as the sampling frequency
decreases, the empirical values get closer to the theoretical limit under
$H_{0}.$ For very high sampling frequencies, the results are consistent with
some mixture of the noise driving the asymptotics. They then slowly settle
down toward the limit corresponding to a null hypothesis of a Brownian
present as the sampling frequency decreases, and the noise presumably becomes
less of a factor.%

\begin{figure}\vspace*{5pt}

\includegraphics{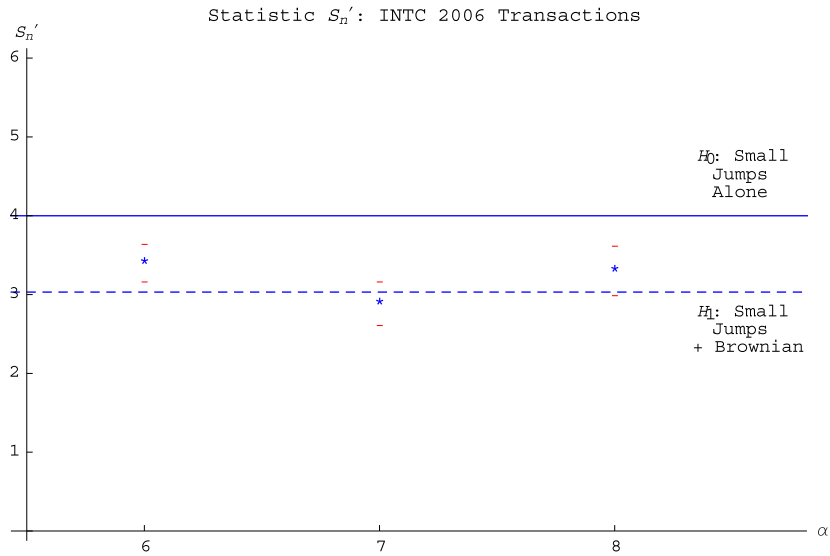}

\caption{Empirical estimates of $S_{n}^{\prime}$ for various truncation levels
$\alpha$ from all Intel transactions during 2006.}\label{fig:S5INTC}
\end{figure}

\begin{figure}

\includegraphics{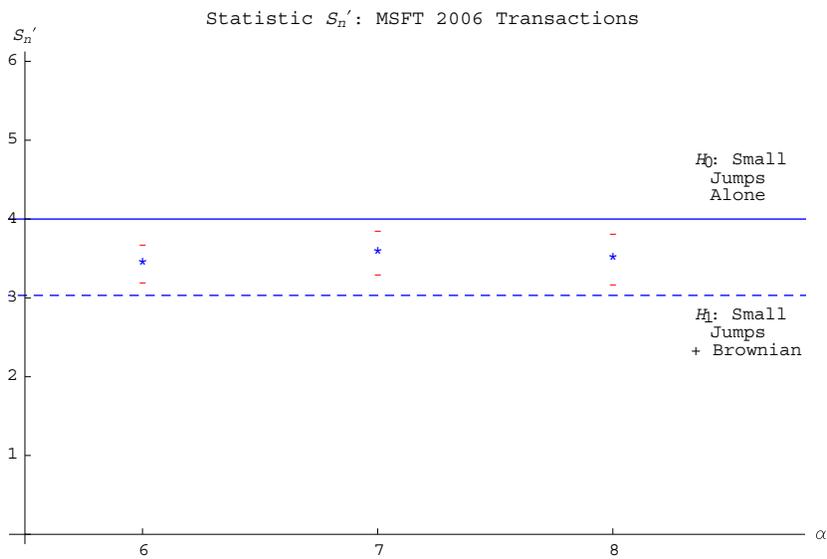}

\caption{Empirical estimates of $S_{n}^{\prime}$ for various truncation levels
$\alpha$ from all Microsoft transactions during 2006.}\label{fig:S5MSFT}\vspace*{5pt}
\end{figure}

Next, we turn to the results of the second test on the same data series in
Figures \ref{fig:S5INTC} and \ref{fig:S5MSFT}. The test statistic
$S_{n}^{\prime}$ is implemented with $\gamma=2,$ with data sampled every
$\Delta_{n}=5$ seconds. The empirical estimates are represented by a star,
with the vertical dashes representing a $95\%$ confidence interval. Also
represented on the plots are the limits corresponding to $H_{0}$ (no Brownian)
and $H_{1}$ (Brownian present). The theoretical limits correspond to those
given in Figures \ref{fig:S5H0Curves} and \ref{fig:S5H1Curves}, except that
the theoretical limit under $H_{1}$\ is plotted for a value of $\beta=1.5$ for
the same reason as above. We find that the empirical estimates tend to be
lower than the value specified by $H_{0},$ which leads to a rejection of the
null hypothesis of no Brownian motion. The estimates are, however, generally
higher than their expected value under $H_{1},$ consistent with the upward
bias identified in simulations, the bias being more pronounced when $\beta$
gets closer to $2$.

To summarize, the answer from both tests appears consistent with the presence
of a continuous component in the data:\ using $S_{n},$ we do not reject the
null of a Brownian motion present, while using $S_{n}^{\prime}$ we reject the
null of its absence.

\section{Technical results}\label{sec:TR}

By a standard localization procedure, we can replace the local
boundedness\ hypotheses in our assumptions by a boundedness assumption, and
also assume that the process $X$ itself, and thus the jump process $\Delta
X_{t}$, are bounded as well. That is, for all results which need Assumption
\ref{ass:A1} we may assume further that, for some constant $C>0$,
\begin{equation}\label{M-1}
|b_{t}|, |\sigma_{t}|, L_{t}, |\Delta X_{t}| \leq C,
\qquad \mbox{hence also }  F_{t}([-C,C]^{c}) = 0 .
\end{equation}
When we need Assumption \ref{ass:A2} we may assume the above, together with
\begin{equation}\label{M-1A}
|\widetilde{b}_{t}|, n_{t}, |\widetilde{\sigma}_{t}|,  \int|x| F_{t}(dx) \leq C.
\end{equation}
We call these \textit{reinforced Assumptions} \ref{ass:A1} or \ref{ass:A2},
and they are assumed in all the sequels instead of mere Assumptions
\ref{ass:A1} or \ref{ass:A2}, according to the case.

Recall that if $\beta<1$, we have (\ref{701}) with $b_{t}^{\prime}$ bounded as
well. Otherwise the decomposition (\ref{701}) is no longer valid, but under reinforced Assumption \ref{ass:A1} we can always write
\begin{equation}\label{702}
Y_{t}^{\prime} = X_{0}+\int_{0}^{t}b_{s}^{\prime\prime} \,ds+\int_{0}^{t}%
\sigma_{s} \,dW_{s}, \qquad  Y_{t}^{\prime\prime} = X_{t}-Y_{t}^{\prime},
\end{equation}
where $b_{t}^{\prime\prime}=b_{t}+\int x1_{\{|x|>1\}} F_{t}(dx)$ defines a
bounded process, and $Y^{\prime\prime}$ is a purely discontinuous martingale.

Also, $K$ below denotes a constant which may change from line to line and may
depend on $C$ above.

The key to all results is clearly the behavior of the processes $B(p,u_{n}%
,\Delta_{n})$ and $U(u_{n},\Delta_{n})$. For establishing this behavior, it is
convenient to introduce a few auxiliary processes, for $u>0$ an arbitrary
cut-off level and $Y$ an arbitrary process
\begin{equation}\label{AU1}
\cases{
\ds B^{\prime}(Y,p,\Delta_{n})_{t} = \sum_{i=1}^{[t/\Delta_{n}]} |\Delta^{n}_{i}
Y|^{p},\cr
\ds D(p,u)_{t} = \sum_{s\leq t}|\Delta X_{s}|^{p} 1_{\{|\Delta X_{s}|\leq
u\}},\cr
\ds D^{\prime}(u)_{t} = \sum_{s\leq t}1_{\{|\Delta X_{s}|>u\}}.
}
\end{equation}

\subsection{Central limit theorems for the auxiliary processes}

This subsection is devoted to recalling or proving some limit theorems for
$B^{\prime}(X,p,\Delta_{n})$ and for the auxiliary processes introduced in
(\ref{AU1}). First, we recall from Theorem 2.4 of \cite{jacod08} that under
Assumption \ref{ass:A1} (and even much more generally),
\begin{equation}\label{ZZ2}
\qquad \cases{
\ds 0<p<2  \quad \Rightarrow\quad  \Delta_{n}^{1-p/2} B^{\prime}(X,p,\Delta_{n})_{t}
\stackrel{\mathbb{P}}{\longrightarrow} A(p)_{t} = m_{p}\int_{0}^{t}|\sigma
_{s}|^{p} \,ds,\cr
p\geq2, X \mbox{ continuous}  \quad \Rightarrow \quad \Delta_{n}^{1-p/2} B^{\prime
}(X,p,\Delta_{n})_{t}  \stackrel{\mathbb{P}}{\longrightarrow}
A(p)_{t},\cr
(u_{n})  \mbox{ satisfies (\ref{M-7})}\quad   \Rightarrow\quad   B(2,u_{n},\Delta
_{n})_{t} \stackrel{\mathbb{P}}{\longrightarrow} A(2)_{t}%
}
\end{equation}
[the last property is proved when $u_{n}=\alpha\Delta_{n}^{\varpi}$ with
$\alpha>0$ and $\varpi\in(0,1/2)$, but the proof works as well when
(\ref{M-7}) holds].

\begin{lemma}
\label{LAU1} Suppose that $X$ is continuous, and let $t\geq0$ and $p>1$ and
$k\geq2$. Under Assumption \ref{ass:A2} the two-dimensional variables
\begin{eqnarray}\label{T-4030}
&&\frac{1}{\sqrt{\Delta_{n}}} \bigl(  \Delta_{n}^{1-p/2}B^{\prime}(X,p,\Delta_{n}
)_{t}-A(p)_{t},\nonumber\\[-8pt]\\[-8pt]
&&\hphantom{\frac{1}{\sqrt{\Delta_{n}}} \bigl(}\Delta_{n}^{1-p/2}B^{\prime}(X,p,k\Delta_{n})_{t}%
-k^{p/2-1}A(p)_{t}  \bigr)  \nonumber
\end{eqnarray}
stably converge in law to a limit which is defined on an extension of
$(\Omega,\mathcal{F},\break(\mathcal{F}_{t})_{t\geq0},\mathbb{P})$ and which,
conditionally on $\mathcal{F}$, is a centered Gaussian variable with
variance--covariance matrix given by
\begin{equation}\label{M-16}
\frac{1}{m_{2p}} \pmatrix{
(m_{2p}-m_{p}^{2})A(2p)_{T} & (m_{k,p}-k^{p/2}m_{p}^{2})A(2p)_{T}\cr
(m_{k,p}-k^{p/2}m_{p}^{2})A(2p)_{T}   & k^{p-1}(m_{2p}-m_{p}^{2})A(2p)_{T}} .
\end{equation}
\end{lemma}

(The same would hold if $p\in(0,1]$, under the additional assumption that
$\sigma_{t}$ is bounded away from $0$.)

\begin{pf*}{Proof of Lemma \ref{LAU1}}
We can assume reinforced Assumption \ref{ass:A2}. The result will follow
from Theorem 7.1 of \cite{jacodsemstat}. Assumption (H) in that paper is
slightly more restrictive than reinforced Assumption \ref{ass:A2}, but a
close look at the proof yields that this theorem still holds in the present situation.

We apply the quoted Theorem 7.1 to the two-dimensional function on
$\mathbb{R}^{k}$ whose components are $|x_{1}|^{p}+\cdots+|x_{k}|^{p}$ and
$|x_{1}+\cdots+x_{k}|^{p}$. This function is $C^{1}$ with derivatives having
polynomial growth. With the notation of that paper, variable
(\ref{T-4030}) with the nontruncated variations is equal to $Z_{n}+R_{n}$,
where
\begin{equation}
Z_{n} = \frac{1}{\sqrt{\Delta_{n}}} \biggl(  \Delta_{n}V^{\prime}(f,k,\Delta
_{n})_{t}-\frac{1}{k}\int_{0}^{t}\rho_{\sigma_{s}}^{\otimes k}(f)\,ds \biggr)
\end{equation}
and $R_{n}$ is a remainder term with second component equal to $0$, and with
first component
\begin{equation}
\Delta_{n}^{1/2-p/2}\sum_{i=k[t/k\Delta_{n}]+1}^{[t/\Delta_{n}]}|\Delta
_{i}^{n}X|^{p}.
\end{equation}
By (\ref{M-1}) we have $E(|\Delta_{i}^{n}X|^{p})\leq K\Delta_{n}^{p/2}$, and
hence, since there are at most $k$ summands in the definition of $R_{n}$, we
deduce that $R_{n}\stackrel{\mathbb{P}}{\longrightarrow}0$. On the other hand,
the aforementioned result yields that $Z_{n}$ converges stably in law to a
limiting variable, which is exactly as described in the statement of the lemma.
\end{pf*}

\begin{lemma}
\label{LAU2} Let $t\geq0$, and suppose Assumption \ref{ass:A1} and $p>\beta$
and $u_{n}\to0$. Then
\begin{equation}
\label{AU4}u_{n}^{\beta-p} D(p,u_{n})_{t} \stackrel{\mathbb{P}}%
{\longrightarrow} \frac{\beta}{p-\beta} \overline{A}_{t},\qquad   u_{n}^{\beta
} D^{\prime}(u_{n})_{t} \stackrel{\mathbb{P}}{\longrightarrow} \overline
{A}_{t}.
\end{equation}
Moreover, if $\beta^{\prime}<\beta/2$ the four-dimensional variables
\begin{equation}\label{AU5}
\pmatrix{
\ds\frac{1}{u_{n}^{\beta/2}} \biggl(  u_{n}^{\beta-p} D(p,u_{n})_{t}- \frac{\beta
}{p-\beta}\overline{A}_{t} \biggr) \cr
\ds\frac{1}{u_{n}^{\beta/2}} \biggl(  (\gamma u_{n})^{\beta-p} D(p,\gamma
u_{n})_{t}-\frac{\beta}{p-\beta}\overline{A}_{t} \biggr) \cr
\ds\frac{1}{u_{n}^{\beta/2}} \bigl(  u_{n}^{\beta} D^{\prime}(u_{n}%
)_{t}-\overline{A}_{t} \bigr) \cr
\ds\frac{1}{u_{n}^{\beta/2}} \bigl(  (\gamma u_{n})^{\beta} D^{\prime}(\gamma
u_{n})_{t}- \overline{A}_{t} \bigr)
}  %
\end{equation}
stably converge in law to a limit which is defined on an extension of
$(\Omega,\mathcal{F} ,\break(\mathcal{F}%
_{t})_{t\geq0},P)$ and which, conditionally on $\mathcal{F}$, is a centered Gaussian variable with variance--covariance matrix
$\overline{A}_{t}\widetilde{C}$, where $\widetilde{C}$ is the $4\times4$
matrix
\begin{equation}\label{AU6}
\widetilde{C}_{t} =   \pmatrix{
\ds\frac{\beta}{2p-\beta} &\ds \frac{\beta \gamma^{\beta-p} }{2p-\beta} & 0 &
0\cr
\ds\frac{\beta \gamma^{\beta-p}}{2p-\beta} & \ds\frac{\beta \gamma^{\beta}}
{2p-\beta} & \ds\frac{\beta(1-\gamma^{\beta-p})}{p-\beta} & 0\cr
0 & \ds\frac{\beta(1-\gamma^{\beta-p})}{p-\beta} & 1 & 1\cr
0 & 0 & 1 & \gamma^{\beta}%
} .
\end{equation}
\end{lemma}

\begin{pf}
Assumption \ref{ass:A1} here implies Assumption 6 of \cite{yacjacod08}, with
the same $\beta$ and $\overline{A}_{t}$, and with $\beta^{\prime}$ there
substituted with any number in $(\beta^{\prime},\beta)$ here. Then all
statements concerning $D(p,u_{n})_{t}$ are in Proposition 5 of that paper.
However, we must redo the proof to obtain the joint convergence for the
processes $D(p,u_{n})$ and $D^{\prime}(u_{n})$.

Let $\widetilde{D}(p,u)$ and $\widetilde{D}^{\prime}(u)$ be the predictable
compensators of $D(p,u)$ and $D^{\prime}(u)$, and $M(u)=u^{\beta
-p}(D(p,u)-\widetilde{D}(p,u))$ and $M^{\prime}(u)=u^{\beta}(D^{\prime
}(u)-\widetilde{D}^{\prime}(u))$. Observe that $\widetilde{D}^{\prime}%
(u)_{t}=\int_{0}^{t} F_{s}([-u,u]^{c})\,ds$ and $ |F_{t}([-v,v]^{c}%
)-v^{-\beta}A_{t} |\leq KL_{t}v^{-\beta^{\prime}}$ by Assumption
\ref{ass:A1}. Therefore, exactly as in the  paper (and (C.23)
and (C.24) in it), we see that if $q>\beta$,
\begin{equation}
\label{AU60}
\cases{
\ds\beta^{\prime}<\beta \quad \Rightarrow\quad   u_{n}^{\beta-q} \widetilde{D}(q,u_{n}
)_{t} \to \frac{\beta \overline{A}_{t}} {q-\beta},\qquad   u_{n}^{\beta
}\widetilde{D}^{\prime}(u_{n})_{t} \to \overline{A}_{t},\cr
\ds\beta^{\prime}<\frac{\beta}2 \quad \Rightarrow\quad  \frac1{u_{n}^{\beta/2}}
 \biggl|u_{n}^{\beta-q} \widetilde{D}(q,u_{n})_{t}-\frac{\beta \overline{A}%
_{t}} {q-\beta} \biggr|\to0, \cr
\hspace*{107.5pt}\ds\frac1{u_{n}^{\beta/2}} | u_{n}^{\beta
}\widetilde{D}^{\prime}(u_{n})_{t}-\overline{A}_{t} |\to0.
}
\end{equation}

The processes $M(u)$ and $M^{\prime}(u)$ are martingales, and if $u\leq v$ the
brackets are given by the following formulas:
\begin{eqnarray*}
\langle M(u),M(v)\rangle &=& (uv)^{\beta-p}\widetilde{D}(2p,u), \qquad \langle
M^{\prime}(u),M^{\prime}(v)\rangle = (uv)^{\beta}\widetilde{D}^{\prime
}(v),\\
\langle M(u),M^{\prime}(v)\rangle &=& 0, \qquad \langle M^{\prime}(u),M(v)\rangle
 = u^{\beta}v^{\beta-p} \bigl(\widetilde{D}(p,v)-\widetilde{D}(p,u) \bigr).
\end{eqnarray*}
This, applied with $(u,v)$ equal to $(u_{n},u_{n})$ or $(u_{n},\gamma u_{n})$
or $(\gamma u_{n},\gamma u_{n})$, and combined with the first part of
(\ref{AU60}), yield that the bracket matrix at time $t$ of the $4$-dimensional
continuous martingale $\overline{M}{}^{n}=u_{n}^{-\beta/2} (M(u_{n}),M(\gamma
u_{n}),\break M^{\prime}(u_{n}),M^{\prime}(\gamma u_{n}))$ converges to $\overline
{A}_{t} \widetilde{C}$ in probability, where $\widetilde{C}$ is given by
(\ref{AU6}). Then as in Proposition~5 of \cite{yacjacod08} one deduces that
$\overline{M}{}_{t}^{n}$ converges stably in law to the limit described in the
statement of the lemma. It remains to deduce from the second part of
(\ref{AU60}) that the difference between $\overline{M}{}_{t}^{n}$ and the
variable defined by (\ref{AU5}) goes to $0$ in probability.
\end{pf}

\subsection{The behavior of $B(p,u_{n},\Delta_{n})_{T}$}

In this subsection we establish the behavior of $B(p,u_{n},\Delta_{n})$ for
the relevant values of $p$ and for the cases not covered by (\ref{ZZ2}). This
is done in several lemmas.

\begin{lemma}
\label{LAU3} Under Assumption \ref{ass:A1}, and if $u_{n}$ satisfies
(\ref{M-7}), we have
\begin{equation}
\label{AU10}B(4,u_{n},\Delta_{n})_{t} \stackrel{\mathbb{P}}{\longrightarrow} 0.
\end{equation}

\end{lemma}

\begin{pf}
We first observe that $B(4,v,\Delta_{n})_{T}$ converges in probability to
$G(v)_{T}=\sum_{s\leq T}|\Delta X_{s}|^{4} 1_{\{|\Delta X_{s}|\leq v\}}$ for
any fixed $v>0$ such that $P(\exists s\leq T\dvtx\break|\Delta X_{s}|=v)=0$. Hence there
is a sequence $v_{m}\rightarrow0$ such that $B(4,v_{m},\Delta_{n})_{T}$
converges in probability to $G(v_{m})_{T}$. On the one hand $B(4,u_{n}%
,\Delta_{n})_{T}\leq\break B(4,v_{m},\Delta_{n})_{T}$ as soon as $u_{n}\leq v_{m}$.
On the other hand we have $G(v_{m})_{T}\rightarrow0$ as $m\rightarrow\infty$.
Then the result follows.
\end{pf}

\begin{lemma}
\label{LAU4} Assume (\ref{M-7}) and reinforced Assumption \ref{ass:A1},
and let $p>0$. If either $p\leq2$, or $p>2$ with $\rho_{-}>\frac
{p-2}{2p-2\beta}$, we have
\begin{equation}
\label{AU11}\Delta_{n}^{1-p/2}  B(p,u_{n},\Delta_{n})_{t}  \stackrel
{\mathbb{P}}{\longrightarrow} A(p)_{t}.
\end{equation}

\end{lemma}

\begin{pf}
We consider decomposition (\ref{702}). In view of (\ref{ZZ2}), it is
enough to prove that under the conditions of the lemma we have
\begin{equation}\label{AU12}
\Delta_{n}^{1-p/2} \bigl(  B(p,u_{n},\Delta_{n})_{t}-B^{\prime}(Y^{\prime
},p,\Delta_{n})_{t} \bigr)   \stackrel{\mathbb{P}}{\longrightarrow} 0.
\end{equation}
The left-hand side above is $\Delta_{n}^{1-p/2}\sum_{i=1}^{[t/\Delta_{n}]}\zeta
_{i}^{n}$, where
\[
\zeta_{i}^{n} = |\Delta_{i}^{n}Y^{\prime}+\Delta_{i}^{n}Y^{\prime\prime}%
|^{p} 1_{\{|\Delta_{i}^{n}X|\leq u_{n}\}}-|\Delta_{i}^{n}Y^{\prime}|^{p}.
\]
With $\kappa=1$ when $p>1$ and $\kappa=0$ otherwise, we have the following
inequalities, for all $m,q>0$:
\begin{eqnarray}\label{T-406}
&&|\Delta_{i}^{n}Y^{\prime}|\geq\frac{u_{n}}{2} \quad  \Rightarrow \quad  |\zeta_{i}%
^{n}|\leq K |\Delta_{i}^{n}Y^{\prime}|^{p+q}/u_{n}^{q},\nonumber\\
&&|\Delta_{i}^{n}X|>2u_{n},\qquad  |\Delta_{i}^{n}Y^{\prime}|\leq\frac{u_{n}}{2}\\
&&\qquad \Rightarrow \quad  |\zeta_{i}^{n}|\leq|\Delta_{i}^{n}Y^{\prime}|^{p}|\Delta
_{i}^{n}Y^{\prime\prime}|^{m}/u_{n}^{m},\nonumber\\
&&|\Delta_{i}^{n}X|\leq2u_{n},\qquad  |\Delta_{i}^{n}Y^{\prime}|\leq\frac{u_{n}%
}{2} \nonumber\\
&&\qquad \Rightarrow \quad  |\zeta_{i}^{n}|\leq K \bigl(  (|\Delta_{i}^{n}%
Y^{\prime\prime}|\wedge u_{n})^{p} \nonumber\\
&  &\hphantom{\qquad \Rightarrow \quad |\zeta_{i}^{n}|\leq} {}   + \kappa|\Delta_{i}^{n}Y^{\prime}|^{p-1}(|\Delta_{i}%
^{n}Y^{\prime\prime}|\wedge u_{n}) \bigr),\nonumber
\end{eqnarray}
where we have used the inequality $ \vert |x+y|^{p}-|x|^{p} \vert
\leq K(|y|^{p}+|x|^{p-1}|y|)$ when $p>1$ and $ \vert |x+y|^{p}%
-|x|^{p} \vert \leq|y|^{p}$ when $p\leq1$. In view of (\ref{M-1}), we
have the estimates
\begin{equation}\label{T-405}
 \cases{
\mathbb{E}(|\Delta_{i}^{n}Y^{\prime\prime}|^{2}) \leq K\Delta_{n},\qquad
q>0 \quad \Rightarrow \quad \mathbb{E}(|\Delta_{i}^{n}Y^{\prime}|^{q}) \leq K_{q}%
\Delta_{n}^{q/2},\cr
r\in(\beta,2]  \quad \Rightarrow \quad  \mathbb{E} \bigl(  (|\Delta_{i}^{n}Y^{\prime\prime
}|\wedge u_{n})^{2} \bigr)   \leq K_{r}\Delta_{n} u_{n}^{2-r}%
}
\end{equation}
(the first estimate is obvious and the second one follows from
Burkholder--Davis--Gundy inequality; the third one follows from (6.25) of
\cite{jacodsemstat} applied to the process $Y^{\prime\prime}$ and with
$\alpha_{n}=u_{n}/\sqrt{\Delta_{n}}$, which goes to $\infty$ by (\ref{M-7}),
and with $r$ as above). Then, using H\"{o}lder's inequality and $(|x|\wedge
u_{n})^{p}\leq u_{n}^{p-2}(|x|\wedge u_{n})^{2}$ when $p>2$, we deduce from
(\ref{T-406}) applied with $q=m=1$ and from $u_{n}\leq K$ that
\[
\Delta_{n}^{1-p/2}\mathbb{E}(|\zeta_{i}^{n}|) \leq  \cases{
\ds K\Delta_{n} \biggl(  \frac{\Delta_{n}^{1/2}}{u_{n}}+u_{n}^{p(1-r/2)}+\kappa
u_{n}^{1-r/2} \biggr),  &\quad if $ p\leq2$,\cr
\ds K\Delta_{n} \biggl(  \frac{\Delta_{n}^{1/2}}{u_{n}}+\Delta_{n}^{1-p/2}%
 u_{n}^{p-r}+u_{n}^{1-r/2} \biggr),     &\quad  if $ p>2$.
}
\]
We have $\Delta_{n}^{1/2}/u_{n}\rightarrow0$ by (\ref{M-7}), hence $E (
\Delta_{n}^{1-p/2}\sum_{i=1}^{[t/\Delta_{n}]}|\zeta_{i}^{n}| )
\rightarrow0$, as soon as $p\leq2$, or $p>2$ and $\rho_{-}\geq\frac
{p-2}{2(p-r)}$. Since $r$ is arbitrary in $(\beta,2]$, we deduce the result.
\end{pf}

\begin{lemma}\label{LAU5}
Let $p\in(0,2]$, and assume reinforced Assumption
\ref{ass:A1} with $\beta<1$ and (\ref{M-7}) with further $\rho_{-}>\frac
{p-1}{2p-2\beta}$ when $p\geq1$. Then, with $X^{\prime}$ given by (\ref{701}),
we have
\begin{equation}
\label{AU13}\Delta_{n}^{1/2-p/2} \bigl(  B(p,u_{n},\Delta_{n})_{t}-B^{\prime
}(X^{\prime},p,\Delta_{n})_{t} \bigr)   \stackrel{\mathbb{P}}{\longrightarrow
} 0.
\end{equation}
\end{lemma}

\begin{pf}
The proof of this lemma is similar to that of the previous one. The left-hand
side of (\ref{AU13}) is $\Delta_{n}^{1/2-p/2}\sum_{i=1}^{[t/\Delta_{n}]}
\zeta_{i}^{n}$, where
\[
\zeta_{i}^{n} = |\Delta_{i}^{n}X^{\prime}+\Delta_{i}^{n}X^{\prime\prime}
|^{p} 1_{\{|\Delta_{i}^{n}X|\leq u_{n}\}}-|\Delta_{i}^{n}X^{\prime}|^{p}.
\]
Then (\ref{T-406}) holds with $(X^{\prime},X^{\prime\prime})$ instead of
$(Y^{\prime},Y^{\prime\prime})$, whereas (\ref{T-405}) is replaced by
\begin{equation}\label{T-407}
 \cases{
\mathbb{E}(|\Delta_{i}^{n}X^{\prime\prime}|) \leq K\Delta_{n},\qquad
q>0 \quad \Rightarrow \quad \mathbb{E}(|\Delta_{i}^{n}X^{\prime}|^{q}) \leq K_{q}
\Delta_{n}^{q/2},\cr
r\in(\beta,1) \quad  \Rightarrow  \quad \mathbb{E}(|\Delta_{i} ^{n}X^{\prime\prime}|\wedge
u_{n}) \leq K_{r} \Delta_{n} u_{n}^{1-r}%
}
\end{equation}
(we now use (6.26) of \cite{jacodsemstat} applied with $\alpha_{n}=u_{n}%
/\sqrt{\Delta_{n}}$ and $r$ as above). Hence, using (\ref{T-406}) for the pair
$(X^{\prime},X^{\prime\prime})$, plus the fact that $(|x|\wedge u_{n})^{p}
\leq u_{n}^{p-m}(|x|\wedge u_{n})^{m}$ for $0<m\leq p$ and H\"{o}lder's
inequality, we deduce that for all $q>0$ and $m\in(0,1)$ and $r\in(\beta,1)$,
and with $\kappa$ as in the previous proof,
\begin{eqnarray*}
&&\Delta_{n}^{1/2-p/2}\mathbb{E}(|\zeta_{i}^{n}|) \\
& &\qquad   \leq K_{r} \Delta
_{n} \biggl(  \frac{\Delta_{n}^{q/2-1/2}}{u_{n}^{q}}+\frac{\Delta_{n}^{m-1/2}%
}{u_{n}^{m}}+\Delta_{n}^{1/2-p/2}u_{n}^{p-r}+\kappa\Delta_{n}^{m-1}%
u_{n}^{1-mr} \biggr) \\
& &\qquad     \leq K_{r}\Delta_{n} (\Delta_{n}^{v_{1}}+\Delta_{n}^{v_{1}}+ \Delta
_{n}^{v_{3}}+\kappa\Delta_{n}^{v_{4}} ),
\end{eqnarray*}
where $v_{1}=q(\frac{1}{2}-\rho_{+})-\frac{1}{2}$ and $v_{2}=m(1-\rho_{+}
r)-\frac{1}{2}$ and $v_{3}=\frac{1-p}{2}+(p-r)\rho_{-}$ and $v_{4}
=m-1+\rho_{-}(1-mr)$. Since $\rho_{+}<1/2$ we have $v_{1}>0$ for $q$ large
enough. When $r\downarrow\beta$ and $m\uparrow1$, we have $v_{2}\to v^{\prime
}_{2}=(1-\rho_{+}\beta)-\frac12$ and $v^{\prime}_{3}\to v^{\prime}_{3}%
=\frac{1-p}2+(p-\beta)\rho_{-}$ and $v_{4}\to v^{\prime}_{4}=(1-\beta)\rho
_{-}$, and (\ref{AU13}) will follow from $v^{\prime}_{j}>0$ for $j=2,3$ and
also for $j=4$ when $p>1$. We have $v^{\prime}_{2}>0$ because $\beta\rho
_{+}<\frac12$. When $p<1$ we have $v^{\prime}_{3}>0$. When $p=1$ then
$v^{\prime}_{3}>0$ if $\rho_{-}>0$, and when $p>1$ we have $v^{\prime}_{3}>0$
and $v^{\prime}_{4}>0$ as soon as $\rho_{-}>\frac{p-1}{2p-2\beta}$. So
(\ref{AU13}) is proved.
\end{pf}

The previous lemma essentially gives the behavior of $B(p,u_{n},\Delta_{n})$
when the leading term is due to the continuous martingale part of $X$. When
this part vanishes, we have another type of behavior, which we describe now.

\begin{lemma}
\label{LAU6} Let $p>1$, and assume reinforced Assumption \ref{ass:A1}.
\begin{longlist}[(ii)]
\item[(i)] If $p>\beta$ and (\ref{M-7}) holds with $\rho_{+}<\frac{p-1}p$ we have
\begin{equation}
\label{T-498}u_{n}^{\beta-p} \bigl(  B(p,u_{n},\Delta_{n})_{t}-D(p,u_{n}%
)_{t} \bigr)  \stackrel{\mathbb{P}}{\longrightarrow} 0
\qquad \mbox{on the set }  \Omega_{t}^{\mathit{noW}}.
\end{equation}

\item[(ii)] If $p\geq2$ and (\ref{M-7}) holds with $\rho_{+}\leq\frac{2-\beta}
{3\beta}$, and if $\beta^{\prime}<\beta/2$, we have
\begin{equation}
\label{T-499}u_{n}^{\beta/2-p} \bigl(  B(p,u_{n},\Delta_{n})_{t} -D(p,u_{n}%
)_{t} \bigr)  \stackrel{\mathbb{P}}{\longrightarrow} 0\qquad
\mbox{on the set }  \Omega_{t}^{\mathit{noW}}.
\end{equation}
\end{longlist}
\end{lemma}

\begin{pf}
Since the variables $B(p,u_{n},\Delta_{n})_{t}$ are the same on the set
$\Omega_{t}^{\mathit{noW}}$ when they are computed on the basis of $X$ or on the basis
of the process $X_{t}-\int_{0}^{t} \sigma_{s} \,dW_{s}$, it is no restriction
to assume that $\sigma_{s}=0$ identically.

The proof is based on the result of \cite{yacjacod08}, when $\sigma_{t}=0$
identically. We have Assumption 7 of that paper with $H=\beta$ and
$a=1-\beta^{\prime}/\beta$ and thus $\phi^{\prime}(x)=x^{-\beta^{\prime}}$. We
can then apply Lemmas 8 of that paper with the version of $\eta(p)_{n}$ given
at the end of Lemma 7 (because $X^{c}=0$ here), to obtain that for
$p>1\vee\beta$ and if $\rho_{+}<\frac{p-1}{p}$ and for any $r\in
 (0,\frac{2}{3\rho_{+}\beta}-\frac{2}{3} )$,
\begin{equation}\label{T-411}
\mathbb{E} \bigl(|B(p,u_{n},\Delta_{n})_{t}-D(p,u_{n})_{t}| \bigr) \leq
 K_{r} t u_{n}^{p-\beta} \eta(p)_{n},
\end{equation}
where
\begin{eqnarray*}
&&\eta(p)_{n} = \sum_{j=1}^{5}(u_{n})^{x_{j}},\\
&&  \cases{
\ds x_{1}=\frac{1}{\rho_{+}}-\beta(1+r),\qquad   x_{2}=\frac{2}{\rho_{+}}-\beta
(2+3r),\cr
\ds x_{3}=r \biggl(1-\frac{\beta}{p} \biggr), \qquad  x_{4}=\frac{p-1}{p\rho_{+}}+\frac{\beta
}{p}-1, \qquad  x_{5}=\beta-\beta^{\prime}.
}
\end{eqnarray*}

Clearly, (\ref{T-498}) follows from (\ref{T-411}), as soon as we can choose
$r\in (0,\frac{2}{3\rho_{+}\beta}-\frac{2}{3} )$ such that $x_{j}>0$ for
all $j=1,\ldots,5$: this is obvious when $\beta<p$ and $\rho_{+}\leq\frac
{p-1}p$.

As for (\ref{T-499}), it will also follow from (\ref{T-411}) if we can choose
$r$ as above, such that $x_{j}>\beta/2$ for all $j=1,\ldots,5$. This property
holds for $x_{5}$ because $\beta^{\prime}<\beta/2$ is assumed, and for $x_{4}$
because $\rho_{+}<1/2$. For $j=1,2,3$, and since $x_{1}$ and $x_{2}$ do not
depend on $p$ and $x_{3}$ increases with $p$, it is enough to consider the
case $p=2$. Then if we let $r$ decrease strictly to $\frac{\beta}{2-\beta}$,
we see that $x_{3}>\beta/2$, whereas $x_{1}$ and $x_{2}$ increase to $\frac
{1}{\rho_{+}}-\frac{2\beta}{2-\beta}$ and to $\frac{2}{\rho_{+}}-\frac
{\beta(4+\beta)}{2-\beta}$ respectively, and these quantities are strictly
bigger than $\beta/2$ if $\rho_{+}$ is strictly smaller than $\frac{4-2\beta}
{\beta(6-\beta)}$ and $\frac{8-4\beta}{\beta(10+\beta)}$. Now, recall that one
should also have $\frac{\beta}{2-\beta}<r<\frac{2}{3\rho_{+}\beta}-\frac{2}%
{3}$, which is possible if and only if $\rho_{+} <\frac{4-2\beta}%
{\beta(4+\beta)}$. All these conditions on $\rho_{+}$ are ensured if $\rho
_{+}\leq\frac{2-\beta}{3\beta}$.
\end{pf}

\subsection{The behavior of $U(u_{n},\Delta_{n})$}

The behavior of $U(u_{n},\Delta_{n})$ has been exhibited in \cite{yacjacod09b},
including a central limit theorem. However, here we need a joint CLT, at
least on the set $\Omega_{T}^{\mathit{noW}}\cap\Omega_{T}^{i\beta}$, for the pair
$(U(u_{n},\Delta_{n})_{T},\break B(2,u_{n},\Delta_{n}))$, and even for this pair
jointly with the similar pair with the truncation levels $\gamma u_{n}$. For
this we will use Lemma \ref{LAU2}, and we thus need to show that the
difference $U(u_{n},\Delta_{n})-D^{\prime}(u_{n})$ is negligible, after a
suitable normalization. To this effect, we use the contorted way of using the
aforementioned CLT for $U(u_{n},\Delta_{n})_{T}$, but knowing this result it
seems the shortest route toward the desired joint CLT.

\begin{lemma}
\label{LAU7} Assume reinforced Assumption \ref{ass:A1}.
\begin{longlist}[(ii)]
\item[(i)] Under (\ref{M-7}) we have
\begin{equation} \label{T-451}
u_{n}^{\beta} U(u_{n},\Delta_{n})_{t}  \stackrel{\mathbb{P}}{\longrightarrow
} \overline{A}_{t}.
\end{equation}

\item[(ii)] If moreover $\beta^{\prime\prime}<\frac{\beta}{2+\beta}$ and
$\beta^{\prime}<\frac{\beta}{2}$ and (\ref{M-7}) holds with $\rho_{+}<\frac
{1}{2+\beta}\wedge\frac{2}{5\beta}$, then
\begin{equation} \label{T-452}
u_{n}^{\beta/2}  \bigl(  U(u_{n},\Delta_{n})_{t}-D^{\prime}(u_{n})_{t} \bigr)
 \stackrel{\mathbb{P}}{\longrightarrow} 0.
\end{equation}
\end{longlist}
\end{lemma}

\begin{pf}
In \cite{yacjacod09b} the truncation level was set as $u_{n}=\alpha\Delta
_{n}^{\varpi}$. However, it is obvious that it works with any truncation level
$u_{n}$ subject to (\ref{M-7}), with the conditions on $\varpi$ replaced by
exactly the same conditions on $\rho_{+}$. With this in view, (i) follows from
Proposition 1 of that paper. The proof of (ii) is much more involved, and
broken into several steps.

\textit{Step} (1) We write $U(u_{n},\Delta_{n})_{t}-D^{\prime}(u_{n})_{t}$ as
$H(1)_{t}^{n}+H(2)_{t}^{n}-H(3)_{t}^{n}$, where $H(3)_{t}^{n}=D^{\prime}%
(u_{n})_{t}-D^{\prime}(u_{n})_{\Delta_{n}[t/\Delta_{n}]}$ and $H(j)_{t}%
^{n}=\sum_{i=1}^{[t/\Delta_{n}]}\zeta(j)_{i}^{n}$ for $j=1,2$, with
\begin{eqnarray*}
\zeta(1)_{i}^{n}&=&1_{\{\Delta_{i}^{n}D^{\prime}(u_{n})=0, |\Delta_{i}%
^{n}X|>u_{n}\}},\\
\zeta(2)_{i}^{n}&=&1_{\{\Delta_{i}^{n}D^{\prime}(u_{n}%
)\geq1, |\Delta_{i}^{n}X|>u_{n}\}}-\Delta_{i}^{n}D^{\prime}(u_{n}).
\end{eqnarray*}
In this step we prove
\begin{equation}
u_{n}^{\beta/2} H(3)_{t}^{n} \stackrel{\mathbb{P}}{\longrightarrow} 0.
\label{ZZ101}%
\end{equation}
The left-hand side above is nonnegative, with expectation $\mathbb{E} (
\widetilde{D}^{\prime}(u_{n})_{t}-\break\widetilde{D}^{\prime}(u_{n})_{\Delta
_{n}[t/\Delta_{n}]} )  $, which is smaller than $K\Delta_{n}/u_{n}%
^{\beta/2+\beta^{\prime}}$ (see the proof of Lemma~\ref{LAU2}). Since
$\rho_{+}(\beta/2+\beta^{\prime})<3\rho_{+}\beta/2<1$ we deduce (\ref{ZZ101}).

\textit{Step} (2) Let us assume for a moment that we have
\begin{equation}
u_{n}^{\beta/2} H(2)_{t}^{n} \stackrel{\mathbb{P}}{\longrightarrow} 0.
\label{ZZ1}%
\end{equation}

In Proposition 2 of \cite{yacjacod09b}, and upon replacing $\alpha\Delta
_{n}^{\varpi}$ by $u_{n}$, it is proved that under our assumptions on
$\beta^{\prime}$, $\beta^{\prime\prime}$ and $\rho_{+}$, the sequence
$Z_{n}=u_{n}^{-\beta/2} (  u_{n}^{\beta} U(u_{n},\Delta_{n}%
)_{t}-\overline{A}_{t} )  $ converges in law to a limiting variable
$\overline{W}_{t}$ which is centered. On the other hand, Lemma \ref{LAU2}
yields that $Z_{n}^{\prime}=u_{n}^{-\beta/2} (  u_{n}^{\beta} D^{\prime
}(u_{n})_{t}-\overline{A}_{t} )  $ converges in law to a limiting
variable $\overline{W}_{t}^{\prime}$ which is also centered (and, indeed, has
the same law as $\overline{W}_{t}$).

Up to taking a subsequence, assume that the pair $(Z_{n},Z_{n}^{\prime})$
converges in law to a pair $(Z,Z^{\prime})$ of variables which are centered,
whereas $Z_{n}-Z_{n}^{\prime}=u_{n}^{\beta/2}(H(1)_{t}^{n}+H(2)_{t}^{n})$. In
view of (\ref{ZZ1}) it follows that $u_{n}^{\beta/2} H(1)_{t}^{n}$ converges
in law to $Z-Z^{\prime}$. Therefore, since by construction $H(1)_{t}^{n}\geq0$
we must have $Z-Z^{\prime}\geq0$. Since $Z-Z^{\prime}$ is centered, we must
have $Z^{\prime}=Z$ a.s. In other words, for any subsequence of $(Z_{n}%
,Z_{n}^{\prime})$ which converges in law, the limit is a.s. $0$, and by a
subsequence principle it follows that the original sequence $Z_{n}
-Z_{n}^{\prime}$ goes to $0$ in law, hence in probability; this obviously
implies (\ref{T-452}).

At this stage, we are left to prove (\ref{ZZ1})\ which will be implied by the
following:
\begin{equation}
\mathbb{E} (  u_{n}^{\beta/2} |\zeta(2)_{i}^{n}| )   \leq \Delta
_{n} v_{n} \label{ZZ6}%
\end{equation}
for a sequence $v_{n}\rightarrow0$.

We recall the property (B.12) of \cite{yacjacod08}: denoting by $R_{1}%
^{n},\ldots,R_{m}^{n},\ldots$ the successive jump times of $D^{\prime}(u_{n})$
occurring after $(i-1)\Delta_{n}$ (with any fixed $i$), we have $P(R_{j}%
^{n}\leq i\Delta_{n})\leq K^{j} \Delta_{n}^{j} u_{n}^{-j\beta}$. This
implies
\[
\mathbb{P}\bigl(\Delta_{i}^{n}D^{\prime}(u_{n})\geq2\bigr) \leq K\Delta_{n}^{2}%
 u_{n}^{-2\beta},\qquad  \mathbb{E} \bigl(  \Delta_{i}^{n}D^{\prime}%
(u_{n}) 1_{\{\Delta_{i}^{n}D^{\prime}(u_{n})\geq2\}} \bigr)   \leq
 K\Delta_{n}^{2} u_{n}^{-2\beta}.
\]
Since $\rho_{+}<2/(3\beta)$ we have $\Delta_{n}/u_{n}^{3\beta/2}\rightarrow0$.
Therefore, for proving (\ref{ZZ6}) it remains to show that
\begin{equation} \label{ZZ7}
u_{n}^{\beta/2} \mathbb{P} \bigl(  \Delta_{i}^{n}D^{\prime}(u_{n}%
)=1, |\Delta_{i}^{n}X|\leq u_{n} \bigr)   \leq \Delta_{n} v_{n}.%
\end{equation}

Set
\[
X^{\prime\prime}(u_{n})_{t} = \sum_{s\leq t}\Delta X_{s} 1_{\{|\Delta
X_{s}|>u_{n}\}},\qquad   X^{\prime}(u_{n}) = X-X^{\prime\prime}(u_{n}).
\]
We have estimate (B.15) of \cite{yacjacod08} again, with $H=\beta$ and
$\phi^{\prime}(x)=x^{-\beta^{\prime}}$. Thus, since on the set $\{\Delta
_{i}^{n}D^{\prime}(u_{n})=1\}$ the process $X^{\prime\prime}(u_{n})$ is
piecewise constant and with a single jump on the interval $\{(i-1)\Delta
_{n},i\Delta_{n}]\}$, and the size of this jump is bigger than $u_{n}$, we
deduce
\begin{equation}\label{ZZ8}
\qquad \mathbb{P} \bigl(  \Delta_{i}^{n}D^{\prime}(u_{n})=1, |\Delta_{i}^{n}%
X^{\prime\prime}(u_{n})|\leq u_{n}(1+w_{n}) \bigr)   \leq K\Delta_{n} (
u_{n}^{-\beta} w_{n}+u_{n}^{-\beta^{\prime}} )  %
\end{equation}
for any choice of the sequence $w_{n}$ decreasing to $0$.

Finally we use estimate (61) of \cite{yacjacod09b} to obtain for all
$q\geq2$
\begin{equation}\label{ZZ9}
\mathbb{P} \bigl(\Delta_{i}^{n}D^{\prime}(u_{n})=1, |\Delta_{i}^{n}X^{\prime
}(u_{n})|>u_{n}w_{n}\bigr)  \leq K\frac{\Delta_{n}^{2}}{w_{n}^{2}%
 u_{n}^{2\beta}}+K_{q}\frac{\Delta_{n}^{q/2}}{w_{n}^{q}u_{n}^{q}}.
\end{equation}
Of course the left-hand side of (\ref{ZZ7}) is smaller than $u_{n}^{\beta/2}$ times
the sum of the left-hand sides of (\ref{ZZ8}) and (\ref{ZZ9}). Therefore, it remains
to prove that we can choose the sequence $w_{n}$ and $q\geq2$ in such a way
that $y_{n}(j)\rightarrow0$ for $j=1,2,3,4$, where
\begin{eqnarray*}
y_{n}(1)&=&u_{n}^{\beta/2-\beta^{\prime}}, \qquad  y_{n}(2)=\frac{w_{n}}{u_{n}%
^{\beta/2}},\\
y_{n}(3)&=&\frac{\Delta_{n}}{u_{n}^{3\beta/2}w_{n}^{2}}%
,\qquad   y_{n}(4)=\frac{\Delta_{n}^{q/2-1}}{u_{n}^{q-\beta/2}w_{n}^{q}}.
\end{eqnarray*}
We have $y_{n}(1)\rightarrow0$ by hypothesis. Upon taking $w_{n}=u_{n}^{r}$
for some $r$, this amounts to showing that one can find $r>0$ and $q\geq2$
such that $r>\frac{\beta}{2}$ and $\frac{1}{\rho_{+}}-2r>\frac{3\beta}{2}$ and
$q-2>(q(2r+2)-\beta)\rho_{+}$. The last condition is satisfied for $q$ large
enough as soon as $2(r+1)\rho_{+}<1$. Then it is easy to see that the choice
of $r$ is possible if and only if $\rho_{+}<\frac{1}{2+\beta}\wedge\frac
{2}{5\beta}$.
\end{pf}

\subsection{Central limit theorems for $B(p,u_{n},\Delta_{n})$ and
$U(u_{n},\Delta_{n})$}

The previous results allow us to derive joint CLTs for the processes
$B(p,u_{n},\Delta_{n})$ and $U(u_{n},\Delta_{n})$, as required for Theorems
\ref{TTCLT1} and \ref{TTCLT2}. For the first of these two theorems, we use the
following proposition which follows from Lemmas \ref{LAU1} and \ref{LAU5}:

\begin{proposition}
\label{P12} Let $p\in(1,2]$ and $t\geq0$ and $k\geq2$. Under Assumption
\ref{ass:A2} and (\ref{M-7}) with $\rho_{-}>\frac{p-1}{2(p-\beta)}$ the
two-dimensional variables
\[
\frac{1}{\sqrt{\Delta_{n}}} \bigl(  \Delta_{n}^{1-p/2}B(p,u_{n},\Delta_{n}
)_{t}-A(p)_{t},\Delta_{n}^{1-p/2}B(p,u_{n},k\Delta_{n})_{t}-k^{p/2-1}
A(p)_{t} \bigr)
\]
stably converge in law to a limit which is defined on an extension of
$(\Omega,\mathcal{F},\break(\mathcal{F}_{t})_{t\geq0},\mathbb{P})$ and which,
conditionally on $\mathcal{F}$, is a centered Gaussian variable with
variance--covariance matrix given by (\ref{M-16}).
\end{proposition}

For the second theorem, we use the following consequence of Lemmas \ref{LAU2},
\ref{LAU6} and~\ref{LAU7}:

\begin{proposition}
\label{P13} Let $t\geq0$ and $\gamma>1$, and suppose Assumption \ref{ass:A1}.
\begin{longlist}[(iii)]
\item[(i)] If $u_{n}\to0$ we have
\begin{equation}
\label{AU41}u_{n}^{\beta} U(u_{n},\Delta_{n})_{t} \stackrel{\mathbb{P}
}{\longrightarrow} \overline{A}_{t}.
\end{equation}

\item[(ii)] If $p>\beta$ and (\ref{M-7}) holds with $\rho_{+}\leq\frac{p-1}p$, we
have
\begin{equation}
\label{AU40}u_{n}^{\beta-p} B(p,u_{n},\Delta_{n})_{t} \stackrel{\mathbb{P}%
}{\longrightarrow} \frac{\beta}{p-\beta} \overline{A}_{t}
\quad \mbox{\rm in restriction to the set $\Omega_{t}^{\mathit{noW}}$}.
\end{equation}

\item[(iii)] If further $\beta^{\prime\prime}<\frac{\beta}{2+\beta}$ and
$\beta^{\prime}<\frac{\beta}2$, and if (\ref{M-7}) holds with $\rho_{+}%
<\frac1{2+\beta}\wedge\frac2{5\beta}\wedge\frac{2-\beta}{3\beta}$, the
four-dimensional variables
\[
\pmatrix{
\ds\frac{1}{u_{n}^{\beta/2}} \biggl(  u_{n}^{\beta-p} B(p,u_{n},\Delta_{n})_{t}-
\frac{\beta}{p-\beta}\overline{A}_{t} \biggr) \cr
\ds\frac{1}{u_{n}^{\beta/2}} \biggl(  (\gamma u_{n})^{\beta-p} B(p,\gamma
u_{n},\Delta_{n})_{t}-\frac{\beta}{p-\beta}\overline{A}_{t} \biggr) \cr
\ds\frac{1}{u_{n}^{\beta/2}} \bigl(  u_{n}^{\beta} U(u_{n},\Delta_{n})_{t}
-\overline{A}_{t} \bigr) \cr
\ds\frac{1}{u_{n}^{\beta/2}} \bigl(  (\gamma u_{n})^{\beta} U(\gamma u_{n}
,\Delta_{n})_{t}-\overline{A}_{t} \bigr)
}
\]
stably converge in law, in restriction to the set $\Omega_{t}^{\mathit{noW}}$, to a
limit which is defined on an extension of $(\Omega,\mathcal{F}
,(\mathcal{F}_{t})_{t\geq0},P)$ and which, conditionally
on $\mathcal{F}$, is a centered Gaussian variable with
variance--covariance matrix $\overline{A}_{t}\widetilde{C}$, with
$\widetilde{C}$ given by (\ref{AU6}).
\end{longlist}
\end{proposition}

\subsection{Proof of the theorems}

It remains to prove the main theorems, for which we can assume the reinforced
assumptions if necessary, without restriction.

First, the consistency results (\ref{T-6}) and (\ref{T-12}) are obvious
consequences of (\ref{AU11}), (\ref{AU40}) and (\ref{AU41}), plus the facts
that $\overline{A}_{T}>0$ on $\Omega^{i\beta}_{T}$ and $A(p)_{T}>0$ on
$\Omega_{T}^{W}$.

Second, in order to prove Theorem \ref{TTCLT1} we use Proposition \ref{P12}
which, upon using the ``delta method,'' shows
that under the stated assumptions the variables $\frac{1}{\sqrt{\Delta_{n}}%
} (S_{n}-k^{p/2-1})$ converge stably in law, in restriction to $\Omega
_{T}^{W}$, to a variable which conditionally on $\mathcal{F}$ is centered Gaussian with variance
\[
V = N(p,k) \frac{A(2p)_{T}}{(A(p)_{T})^{2}}.
\]
With $V_{n}$ given by (\ref{T-8}), we have $\frac{1}{\Delta_{n}}%
 V_{n}\stackrel{\mathbb{P}}{\longrightarrow}V$ by (\ref{AU11}), and the result
readily follows.

In the same way, Proposition \ref{P13} yields that $\frac{1}{u_{n}^{\beta/2}%
} (S_{n}^{\prime}-\gamma^{2})$ converges stably in law, in restriction to
$\Omega_{T}^{\mathit{noW}}\cap\Omega_{T}^{i\beta}$, to a variable which conditionally
on $\mathcal{F}$ is centered Gaussian with variance
\[
V^{\prime} = \frac{\gamma^{4}}{\overline{A}_{T}}   \biggl(  \frac{\beta
(2-\beta)^{2}}{4-\beta}+1 \biggr)  (1+\gamma^{\beta}-2\gamma^{\beta-2}).
\]
If $V_{n}^{\prime}$ is given by (\ref{CLT V'}), then $\frac{1}{u_{n}^{\beta}%
} V_{n}^{\prime} \stackrel{\mathbb{P}}{\longrightarrow}V^{\prime}$ in
restriction to $\Omega_{T}^{\mathit{noW}}\cap\Omega_{T}^{i\beta}$ by (\ref{AU40}) and
(\ref{AU41}). This finishes the proof of Theorem \ref{TTCLT2}.

Finally, for both Theorems \ref{TT2} and \ref{TT4}, the claims concerning the
asymptotic level of the tests are trivial consequences of two central
limit Theorems \ref{TTCLT1} and \ref{TTCLT2}. It remains to prove that the
asymptotic power is $1$ in both cases. By virtue of (\ref{T-6}) and
(\ref{T-12}), this will follow from the next two properties, under the
appropriate assumptions
\begin{equation}  \label{M223}
 \cases{
V_{n} \stackrel{\mathbb{P}}{\longrightarrow} 0,&\quad  on the set $  \Omega
_{T}^{\mathit{noW}}\cap\Omega_{T}^{i\beta}$,\cr
V_{n}^{\prime} \stackrel{\mathbb{P}}{\longrightarrow} 0,&\quad
on the set $  \Omega_{T}^{W}\cap\Omega_{T}^{i\beta}$.
}
\end{equation}
The first of these properties follows from (\ref{AU40}), and the second one
follows from (\ref{ZZ2}), (\ref{AU10}) and (\ref{AU41}).

\section*{Acknowledgments}

We are very grateful to a referee and an Associate Editor for many helpful
comments.

\printaddresses

\end{document}